\documentclass[12pt]{article}
 
\usepackage{amssymb,amsmath}
\usepackage{graphicx}
\usepackage{epsfig}
\usepackage{color}
\usepackage{cite}
\usepackage[british]{babel}
\include{psfig}

\allowdisplaybreaks

\title{Limit theorems for random spatial drainage networks}

\author{Mathew D.~Penrose\thanks{Partially
 supported by the Alexander von Humboldt Foundation 
through a Friedrich Wilhelm Bessel Research Award} \footnote{e-mail:
 \texttt{m.d.penrose@bath.ac.uk}}
 \\
 \normalsize
 Department of Mathematical Sciences,
 University of Bath,\\
\normalsize
 Claverton Down, Bath BA2 7AY, England.
\and Andrew R.~Wade\thanks{Partially
 supported by the 
 Heilbronn Institute for Mathematical Research.} \footnote{e-mail: \texttt{Andrew.Wade@bris.ac.uk}}
\\
\normalsize
 Department of Mathematics,
 University of Bristol,\\
\normalsize
 University Walk, Bristol BS8 1TW, England.}

\newcommand{\ud}{\mathrm{d}}

\newtheorem{theorem}{Theorem}[section]

\newtheorem{lemma}{Lemma}[section]

\newtheorem{definition}{Definition}[section]
 
\newcommand{\tod}{\stackrel{d}{\longrightarrow}}
\newcommand{\eqd}{\stackrel{d}{=}}
\newcommand{\toP}{\stackrel{P}{\longrightarrow}}

\newcommand{\rems}{\noindent \textbf{Remarks. }}
\newcommand{\proof}{\noindent \textbf{Proof. }}

\def\Exp{{\mathbb{E}}}
\def\Pr{{\mathbb{P}}}
\def\Var{{\mathbb{V}\mathrm{ar}}}
\def\R{\mathbb{R}}
\def\Z{\mathbb{Z}}
\def\1{{\bf 1 }}
\def\bx{{\bf x }}
\def\bz{{\bf z }}
\def\bw{{\bf w }}
\def\by{{\bf y }}
\def\re{{\rm e }}
\def\bn{{\bf n}}

\def\N{\mathbb{N}}
\def\Lm{\mathcal{L}_{\rm max}}
\def\Om{\mathcal{O}_{\rm max}}
\def\NN{\mathcal{N}}
\def\Y{\mathcal{Y}}
\def\po{\preccurlyeq}
\def\pos{\preccurlyeq_*}

\def\0{{\bf 0}}

\def\U{\mathcal{U}}
\def\X{\mathcal{X}}
\def\O{\mathcal{O}}
\def\W{\mathcal{W}}

\def\V{{\cal V}}

\def\tO{\tilde {\cal O}}
\def\bV{\mathbf{V}}
\def\bU{\mathbf{U}}
\def\bY{\mathbf{Y}}

\def\Po{\mathcal{P}}
\def\H{\mathcal{H}}
\def\LL{\mathcal{L}}
\def\card{{\rm  card}}
\def\eps{{\varepsilon}}

 \def\Wal{W_\alpha}

\def\tLalph{\tilde \LL^{d,\alpha}}

\def\tJo{\tilde J_1}                                                                               
\def\tHo{\tilde H_1}
\def\tGo{\tilde G_1}
\def\tJal{\tilde J_\alpha}                                                                               
\def\tHal{\tilde H_\alpha}
\def\tGal{\tilde G_\alpha}
\def\be{{\bf e}}
\def\bm{{\bf m}}

\begin{document}

\maketitle 

\begin{abstract}
Suppose that under the action of gravity,
liquid drains 
 through the
unit $d$-cube via a minimal-length network of channels
constrained to pass through   
random sites and to flow
with nonnegative component
in
one of the canonical orthogonal basis directions of $\R^d$, $d \geq 2$.
The resulting network is a 
version
of the so-called minimal directed spanning tree. 
 We give laws of large numbers and convergence
in distribution results on the large-sample asymptotic behaviour of the
 total power-weighted edge-length of the
network on uniform random points in $(0,1)^d$. 
The distributional results exhibit a 
weight-dependent phase transition between Gaussian 
and boundary-effect-derived distributions.
These boundary
contributions are characterized
in terms of limits of the so-called on-line nearest-neighbour graph,
a natural model of spatial network evolution, for which we also 
present some new results. Also, we
give a convergence in distribution result for the length of the longest edge in the drainage network; 
when $d=2$, the
limit is expressed in terms of Dickman-type variables.
\end{abstract}

\vskip 2mm


\vskip 2mm


\section{Introduction}

We consider a continuum model of drainage through a 
porous medium in $\R^d$ ($d \in \N := \{1,2,3,\ldots\}$), which we first describe
informally. Let $\{\be_1,\ldots,\be_d\}$ be the canonical orthonormal
basis of $\R^d$. We distinguish the
$\be_d$ direction and suppose that `gravity' acts in direction $-\be_d$;
in free space, liquid would fall in exactly the  $-\be_d$ direction.

Informally,
consider a unit $d$-cube, representing a block of 
 porous material. We scatter a certain
finite set $\X$ 
of points in this cube, representing special sites  in the medium. 
We constrain liquid to drain in channels that visit every site
and travel in straight lines from site
to site. The vectors of each channel must have a non-positive component
in the $\be_d$ direction; that is, they must respect gravity. 
The collection of channels spanning $\X$ satisfying these
conditions we call a {\em drainage network} on $\X$.
A natural question is to find the most efficient arrangement of channels
 satisfying the above constraints, i.e., a drainage network that
 is in some sense optimal. As we shall see,
 an answer to this question
 is a version of the so-called {\em
minimal directed spanning tree} (MDST for short) on the vertices $\X$. 

More mathematically, let $\X$ be a finite point set in $(0,1)^d$ whose
points have distinct $d$-th coordinates. 
We construct a directed graph on vertex set $\X$ as follows.
Join each vertex $\bx \in \X$  by a
directed edge to
a Euclidean nearest neighbour (if one exists, and arbitrarily breaking any ties)
 amongst those points $\by \in \X \setminus
\{ \bx\}$ such that $\by \pos \bx$. Here $\pos$ is the order
on $\X$  induced by the order on $d$-th coordinates:
$(x_1, \ldots, x_d) \pos (y_1,\ldots,y_d)$ if and only if $x_d \leq y_d$.
We call the directed graph so constructed
the MDST on $(\X;\pos)$:
it is a mathematical solution
to the problem of constructing a minimal-length
drainage network on $\X$ as
informally described above.

The subject of this paper
is the MDST on $(\Po_n ; \pos)$
where  $\Po_n$ is a homogeneous
Poisson point process of
intensity $n>0$ on
 $(0,1)^d$.
Then (with probability $1$),
$\Po_n$ is indeed a finite point set with
 distinct
$d$-th coordinates so that
the MDST is almost surely well-defined. 
We study the total power-weighted edge-length of the MDST on $(\Po_n;\pos)$
as $n \to \infty$, and also
 the length of the longest edge.

The MDST on $(\Po_n;\pos)$ is an example of a random spatial graph,
that is, a graph 
generated by scattering
points randomly into a region of space and connecting 
them according to some prescribed rule.
Motivated in part by
 real-world networks with
spatial content, such as communications networks (including the Internet),
social networks, and physical networks,
 a substantial body of recent research has dealt with the large-sample
asymptotic theory of such graphs.
 Examples
include the geometric graph, the nearest-neighbour graph,
and the minimal-length spanning tree (MST). 
See, for example,
\cite{avbert,kl,penbook,mdp,mpbern,mpejp,penyuk1,penyuk2,steelebook,yukbook}.
A feature that   distinguishes the MDST considered here
from other random spatial graphs
is that the
 constraint on direction of the edges can lead to significant (indeed, sometimes dominating)
 boundary effects due
to the possibility of long edges occurring near 
the lower boundary cube $(0,1)^{d-1}$ orthogonal to $\be_d$.
Another difference  
is the fact that there is no
uniform upper bound on vertex degrees in the MDST.

In general, the MDST can be defined on any finite partially ordered 
set in $\R^d$, as described in \cite{rooted};
a survey of results on the random MDST
is given in \cite{survey}.
Examples considered previously are
the `cooridnate-wise' (or `South-West') partial ordering
on point sets in $(0,1)^2$ \cite{bhattroy2002,total,rooted}
or in $(0,1)^d$ \cite{blp}, and  the radial spanning tree \cite{bacc}
on  point sets in $\R^2$.
Also, laws of large numbers for the MDST on a class of partial orders of $\R^2$
were given in \cite{llns}.

 In this paper we are concerned with the
`South' partial order $\pos$, which is even a total order,  
on point sets in $\R^d$ with distinct $d$-coordinates.
Our main results, Theorems
\ref{llnthm} and \ref{mainth}, give
laws of large numbers, convergence of expectation,
and distributional
convergence results for  the total
 power-weighted edge-length of the
MDST on $(\Po_n;\pos)$ for $d \geq 2$. We also give a convergence result
for the maximum edge-length in the MDST (Theorem \ref{dickman}).
Our main distributional limit result, Theorem \ref{mainth}, reveals two
regimes of limit behaviour for the total power-weighted edge-length
depending on the power-weighting, in which  the limit law
is either purely normal or given in terms of boundary effects 
characterized as distributional limits
of certain {\em on-line nearest-neighbour graphs}. At a critical point between these two regimes, 
there is a phase transition
at which both effects contribute significantly to the limit law. In order to understand
the boundary effects in the MDST, and its
 longest edge, we 
make use of the fact that near to the boundary, the MDST is well-approximated
by a certain on-line nearest-neighbour graph.

In the on-line nearest-neighbour graph (ONG), each point after the first
in a sequence of
points arriving sequentially
in $\R^d$ is
joined to its nearest
neighbour amongst
those points already present.
The ONG itself is of separate interest as a simple growth model
for random networks, such as the world wide web graph (see \cite{boll}).
The total power-weighted length of the ONG has been studied in
\cite{mdp, ong, llns, ong2}. In the present paper, the ONG arises
as a natural tool for studying the structure of the MDST near to the boundary;
we also prove a new result (Theorem \ref{ongmax})
on the length of the longest edge in the ONG on uniform random
points in $(0,1)^d$.

In the particular case of the total weight of
 the MDST on $(\Po_n ; \pos)$ when $d=2$, which is
 one of the most natural cases,
the boundary contributions to the total
power-weighted edge-length limit laws
 can be characterized in terms 
 the limiting distribution of the total weight of the one-dimensional ONG
  (centred as
 necessary).
Results from \cite{ong} say that  
such a
distribution
 is characterized by a distributional fixed-point equation. Such
  fixed-point equations,
  and the `divide and conquer' algorithms
from which they often arise, are also a subject of considerable recent interest; see, for
example,~\cite{ald,neinrusch,rosler0}.

Mathematically, much of the motivating interest
comes from the desire to further understand the
interplay between stochastic geometry and
distributional fixed points previously more commonly seen in
the analysis of algorithms (see e.g.~\cite{neinrusch}). 
This relationship was first seen in our previous work
\cite{total,llns}
on  limit theorems for the
 length of the `South-West'
MDST in the unit square. The present work adds to this by
 considering the `South' MDST,
for which the fixed-point distributions
which arise are different.
We remain some way from having  
a full description of the limits for all possible partial orders,
other shapes of domain and non-uniform densities.

In \cite{total,llns}, only the case $d=2$ of the `South-West' MDST
was studied. In the present paper
we deal also with higher dimensions.
With fairly straightforward modifications, 
the method used in \cite{total} could be
adapted to prove the $d=2$ case
of our Theorem \ref{mainth} below. 
However, at several points the proofs used in \cite{total}
are not easily adapted to higher dimensions,
and thus we have adopted different proofs;
sometimes these improve or extend ideas
from \cite{total} and sometimes we use 
entirely different techniques.
Another difference is that \cite{total,llns}
made use of general results of Penrose and Yukich
\cite{penyuk1,penyuk2} while in the present paper
we instead use the results of Penrose \cite{mpbern,mpejp}
(see also \cite{mdp05})
which are in several ways more convenient  
for the current application. Thus the
results of the present paper
are of a similar (albeit general-dimensional)
flavour to those
in \cite{total,llns}, but the proofs are different.

Before describing our results in detail, we return to the question
 of motivation. 
General motivation for the MDST is as a model
for a constrained optimal transport network (see e.g.~\cite{survey}).
As has been mentioned elsewhere (e.g.~\cite{bhattroy2002}), the MDST can be
motivated by communications networks. However, in the present case the primary 
 motivation is 
from drainage networks. From  this point of view, our choice of
`South' partial ordering seems the most natural, and
the two most natural choices of $d$ are
 $d=2$ and $d=3$. 
 For further references on the
mathematical modelling of drainage networks, and a related infinite
lattice version of this model, for which rather different properties were studied,
see \cite{grs}; for background on modelling of drainage networks in general, see also \cite{rodriguez}.

\section{Statement of results} 
\label{con}

In this section we give formal definitions of our model
and state our main results. 
Let $d \in \N$.
Let $\X$ be a finite subset of $\R^d$
 endowed with the binary relation
$\pos$, for which $(x_1, \ldots, x_d) \pos (y_1,\ldots,y_d)$ if and only if $x_d \leq y_d$.
Assume that all the elements of $\X$
 have distinct $x_d$-coordinates. Under this assumption, $\pos$
is a partial order on $\X$ (in fact, a total order), and so the MDST that we shall
construct fits
into the theory of the MDST on partially ordered sets given in \cite{rooted,survey}.
Let $\card(\X)$ denote the cardinality (number of elements)
of the set $\X$.

 A {\em minimal element}, or {\em sink}, is a vertex 
$\bx \in \X$ for which there
exists no $\by \in \X \setminus \{\bx\}$ such that $\by \pos \bx$. Thus under our definition
of $\pos$ and our assumption on $\X$, there is a unique sink having strictly minimal
$x_d$-coordinate and which we shall denote $\bm(\X)$. 

For a vertex $\bx \in \X \setminus \{ \bm(\X) \}$, 
we say that $\by \in \X \setminus \{\bx\}$
 is a \emph{directed nearest neighbour} (in the $\pos$-sense)
 of $\bx$ with respect to $\X$
if $\by \pos \bx$ and \[
\|\by-\bx\|_d = \min_{\bz \in \X \setminus \{\bx\} : \bz \pos \bx} \|\bz-\bx\|_d; \]
here and subsequently $\|\cdot\|_d$ denotes the Euclidean norm
on $\R^d$.
For each $\bx \in \X \setminus \{ \bm(\X) \}$ 
let $\bn_\bx := \bn( \bx ; \X)$ denote  a directed nearest neighbour of $\bx$
with respect to $\X$, chosen
arbitrarily if $\bx$ has more than one directed nearest neighbour.
A  minimal directed spanning tree (MDST)
on $(\X;\pos)$, or simply `on $\X$' from now on,
 is a directed graph with vertex set 
$\X$ and edge set $\{ (\bx, \bn_\bx) : \bx \in \X \setminus \{ \bm(\X) \}\}$. That is,
there is an edge from each point other than the sink to a directed nearest neighbour. Hence,
ignoring the directedness of the edges, an MDST on $\X$ 
 is a tree rooted at the sink $\bm(\X)$. Note that
an MDST is also a solution to a 
global optimization problem (see \cite{bhattroy2002,rooted}) --- that is, find a minimal-length
spanning
tree (ignoring directedness of the edges) such that each vertex is 
connected to the sink by a unique
directed path, where directed edges must respect $\pos$.

\begin{figure}[h!]
\begin{center}
\includegraphics[angle=0, width=0.4\textwidth]{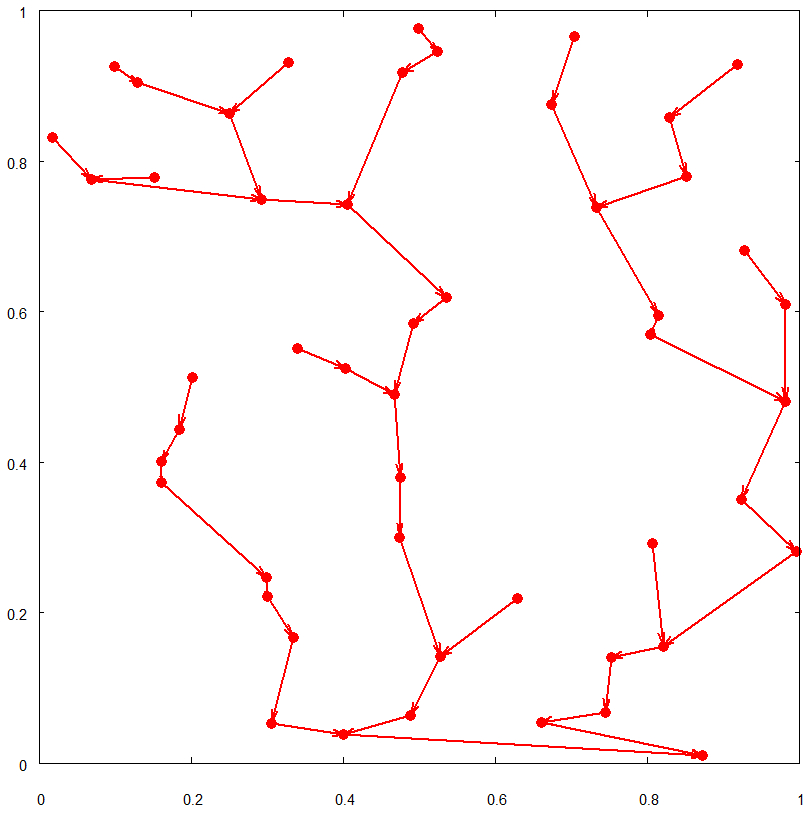}
\includegraphics[angle=0, width=0.5\textwidth]{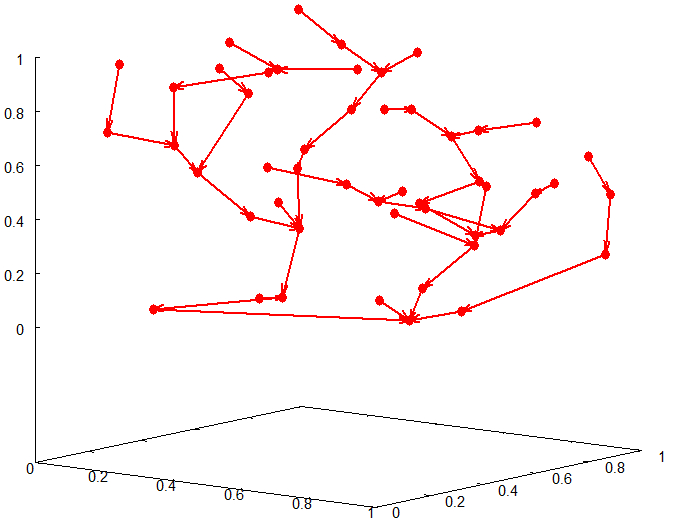}
\end{center}
\caption{Realizations of the MDST under $\pos$
 on 50 simulated uniform random points in
  $(0,1)^2$ (left) and $(0,1)^3$ (right).}
\label{onngfig}
\end{figure}

For $\X \subset \R^d$ with $\card(\X) \geq 2$,
let $d_*(\bx;\X)$ denote the Euclidean distance
from a non-minimal $\bx \in \X$ to a
 directed nearest neighbour $\bn(\bx;\X)$ 
under $\pos$
and   set $d_*(\bm(\X);\X)=0$. 
For $d \in \N$ and $\alpha >0$, define the total
  power-weighted edge-length
of 
 the MDST on $\X$ by
\begin{align*}
 \LL^{d,\alpha} ( \X ) := \sum_{\bx \in \X} (d_*(\bx;\X))^\alpha
 = \sum_{\bx \in \X \setminus \{ \bm(\X) \}} \| \bx - \bn(\bx ; \X) \|_d^\alpha,\end{align*}
 where an empty sum is $0$.
In particular, $\LL^{d,1} (\X)$ is the total Euclidean length of the MDST on $\X$.
Also, define the centred version
$\tilde \LL^{d,\alpha} ( \X ) := \LL^{d,\alpha} ( \X
) - \Exp [ \LL^{d,\alpha} ( \X )
]$. 

From now on we will take $\X$ to be a {\em random} point set
in $(0,1)^d$. In particular, we will take
a homogeneous Poisson point process $\Po_n$ of
intensity $n$ on $(0,1)^d$.
 Note that in
this random setting,
each point of $\Po_n$ almost surely
 has a unique $x_d$-coordinate
 and at most one directed nearest neighbour under $\pos$, so
that $\Po_n$ has a unique MDST, which is rooted at $\bm(\Po_n)$.

We state
and prove
all of our main results in the present paper for the Poisson process $\Po_n$. In all cases,
the authors believe that
analogous results hold for the binomial point process consisting
of $n$ independent uniform random points on $(0,1)^d$
 instead;
it should be possible to use
standard de-Poissonization arguments (such as applied
in similar circumstances in \cite{rooted, total}) to verify this.

In the present paper we are concerned with $d \geq 2$. 
When $d=1$, $\pos$ coincides with the coordinatewise partial order
$\po^*$
(and indeed the total order $\leq$ on $\R$)
and so our `South' MDST is the same as the `South-West'
MDST here.
Moreover, $\LL^{1,\alpha} (\Po_n)$ is a sum
of powers of spacings of 
 uniform points,
and it can be studied using standard
Dirichlet
spacings results (see e.g.~\cite{bill,darling}).
For instance, Darling (see \cite{darling}, p.~245)
essentially 
gives a central limit theorem for the
binomial point process analogue of $\LL^{1,\alpha} (\Po_n)$.
From now on we fix $d \in \{2,3,\ldots\}$.

Our first result describes the first-order behaviour of $\LL^{d,\alpha} (\Po_n)$
as $n \to \infty$. In particular, we have a law of large numbers for 
$\alpha\in(0,d)$, and also asymptotic
results for the expectation when $\alpha \geq d$. 
In $d=2$, the binomial point process
version
 of Theorem \ref{llnthm}(i) is contained in 
 the $\phi = \pi$ case of
 Theorem 5
 of \cite{llns}.
For $d \in \N$, let
\begin{align}
\label{0818c}
v_d
: = \pi^{d/2} \left[ \Gamma \left( 1+ (d/2) \right) \right]^{-1},
\end{align}
 the volume of the unit
$d$-ball (see e.g.~\cite{huang} equation (6.50)); here $\Gamma(\cdot)$ denotes the Euler
Gamma function. 

\begin{theorem} \label{llnthm} 
Suppose $d \in \{2,3,4,\ldots\}$. \begin{itemize}
\item[(i)]
Suppose $\alpha \in (0,d)$.
Then
as $n \to \infty$,
\begin{align}
n^{(\alpha/d)-1} \LL^{d,\alpha} ( \Po_n ) \to
2 ^{\alpha/d} \Gamma (1+(\alpha/d)) v_d^{-\alpha/d} , \textrm{ in } L^1.
\label{0728e}
\end{align}
\item[(ii)] Suppose $\alpha \geq d$. Then there exists
$\mu'(d,\alpha) \in (0,\infty)$ such that,
as $n \to \infty$
\begin{align}
\label{011b}
 \Exp [ \LL^{d,\alpha} (\Po_n)] \to \mu'(d,\alpha) .\end{align}
 Moreover,
 we can express 
 \[ \mu'(d,\alpha) = \mu(d-1,\alpha)  +\1_{\{\alpha =d\}} 2 v_d^{-1} ,\]
 where constants $\mu(d-1,\alpha) \in (0,\infty)$ 
can be characterized in terms of limits of certain on-line
nearest-neighbour graphs:
$\mu(\cdot,\cdot)$ is as given in Proposition 2.1 of
\cite{ong2}; see (\ref{969}) below.
 In particular, for $\alpha \geq 2$
 \[ \mu(1,\alpha) = \frac{2}{\alpha(\alpha+1)} \left( 1+\frac{2^{-\alpha}}{\alpha-1} \right) .\]
\end{itemize}
\end{theorem} 

 One can generalize the statement of 
Theorem \ref{llnthm}(i) to more general point processes under certain conditions; see
\cite{mpbern,mdp05} for a general framework.

Our second main result (Theorem \ref{mainth}, below)
presents convergence in distribution
results for $\LL^{d,\alpha} (\Po_n)$;
the distributional limits contain Gaussian random variables
and also random variables defined as distributional limits
of certain on-line nearest-neighbour graphs (see Section \ref{secong}). In general
we do not give an explicit description of the latter
distributions. However, in
the case of $d=2$, the limits in question can be characterized
as solutions to distributional fixed-point equations, which we describe at the end
of this section.

We now state our main convergence in distribution result.
Let $\NN(0,s^2)$ denote the normal distribution 
with mean zero and variance $s^2 \geq 0$;
included is the degenerate case $\NN(0,0)$.
 By `$\tod$' we denote
convergence in distribution. 

\begin{theorem} \label{mainth}
Suppose $d\in\{2,3,4,\ldots\}$ and $\alpha>0$.
Then
there exists a constant $s_\alpha^2 \in [0,\infty)$
such that,
for a normal random variable 
$\Wal
\sim \NN(0,s_\alpha^2)$,
as $n
\to \infty$:
\begin{align*}
n^{(\alpha/d)-(1/2)} \tilde {\cal L}^{d,\alpha} ( \Po_n )
\tod \Wal ~~~ (0 < \alpha < d/2); 
 \\
\tilde {\cal L}^{d,d/2} ( \Po_n ) \tod W_{d/2} + Q(d-1,d/2) ;~~
 \\
\tilde {\cal L}^{d,\alpha} ( \Po_n ) \tod Q(d-1,\alpha)  ~~~ (\alpha > d/2)
 .
\end{align*}
 Here  
the  $Q(d-1,\alpha)$ 
are mean-zero random variables as given in Lemma \ref{onngthm} below
and independent  of the $W_\alpha$; in particular
$Q(1,\alpha)=\tGal$ for $\alpha \geq 1$, where $\tGal$
has the distribution given by (\ref{0926a}) or (\ref{0926p}) below.
\end{theorem}
\rems 
(a) It can be shown that the limiting variance $s_\alpha^2$
of the normal component in the above limits
is strictly positive for $\alpha>0$, using, for example, techniques similar to those in \cite{avbert} or
\cite{penyuk2} (see
Lemma 6.2 of the extended version of \cite{total}
 for an example of such a result for
 a different MDST model).\\
(b)
The normal random variables $\Wal$
arise from the edges
away from the lower boundary of the $d$-cube (see Section \ref{ltot}).
The variables $Q(d-1,\alpha)$
arise from the edges very close to the boundary, where the MDST
is asymptotically close to a $(d-1)$-dimensional on-line nearest-neighbour graph: this is formalized in
 Section \ref{bdry} below.\\
(c)
 Theorem \ref{mainth} 
 indicates a phase transition
in the character of the limit law as $\alpha$ increases.
The normal contribution  
dominates for $\alpha \in (0,d/2)$, while the boundary contribution dominates
for $\alpha > d/2$. In the critical case $\alpha=d/2$ (such as the natural
case $d=2$ and $\alpha=1$) neither effect dominates
and both terms contribute significantly to the asymptotic behaviour. The intuition here is that
increasing $\alpha$ increases the relative importance
of long edges, such as, typically, those
near to the boundary. \\
(d) As will be demonstrated below (see Lemma \ref{onngthm}), the random variables $Q(d-1,\alpha)$
can be characterized as distributional limits of certain on-line nearest-neighbour graphs. 
It is known (see \cite{ong}) that the $Q(d-1,\alpha)$
are non-Gaussian for $\alpha>d-1$. When $d=2$ much more is known (see \cite{ong});
$Q(1,\alpha)$ can be characterized in terms of a distributional fixed-point equation (see (\ref{0926a})
and (\ref{0926p}) below). In particular, $Q(1,\alpha)$ is non-Gaussian for $\alpha \geq 1$. The authors suspect
that for general $d$, $Q(d-1,\alpha)$ is in fact non-Gaussian for all $\alpha \geq d/2$. \\

Theorem \ref{dickman} below gives a convergence in distribution
 result on the length of the longest edge in the MDST  
on 
$(\Po_n;\pos)$. A similar result (in $d=2$ only)
 for the longest edge in the `South-West' MDST
was given in \cite{rooted}.
Let $\Lm^d (\X)$ denote the length of the longest edge in the MDST (under $\pos$) 
on point set $\X \subset (0,1)^d$:
\[ \Lm^d (\X) := \max_{\bx \in \X} d_* ( \bx ; \X) = 
\max_{\bx \in \X \setminus \{ \bm(\X) \}} \| \bx - \bn( \bx; \X) \|_d .\]
In the particular case $d=2$,
the distributional limit arising in Theorem \ref{dickman} below
is expressed in terms of the {\em max-Dickman} distribution (named
after Dickman's work \cite{dickman}
on the asymptotic distribution of large prime factors), which can be characterized
as the distribution of a random variable $M$ satisfying the fixed-point equation
\begin{align}
\label{maxdick}
 M \eqd \max \{ 1-U, UM \},\end{align}
where $U$ is uniform on $(0,1)$
and independent of the $M$ on the right. (Here and subsequently `$\eqd$' 
denotes equality in distribution.)
See \cite{rooted,survey} and references therein
for more information on the max-Dickman distribution; it has appeared in many contexts, and
a picture of part of its density function is on the front cover of the second edition of
 Billingsley's book \cite{bill}. In particular, we note that
 $M$ can be characterized as the first component of the Poisson--Dirichlet distribution with
 parameter 1, and $\Exp[M] \approx 0.6243299$ is 
 Dickman's constant 
 (see 
 \cite{dickman} p.~9).

\begin{theorem}
\label{dickman}
Let $d \in \{2,3,\ldots,\}$.
There exists a random variable $Q_{\rm max}(d-1)$
such that
\begin{align*}
\Lm^d (\Po_n)  \tod
Q_{\rm max}(d-1),\end{align*}
as $n \to \infty$. Moreover,
$Q_{\rm max}(d-1)$ is characterized in terms of the ONG (see Theorem
\ref{ongmax} below); in particular
\[ Q_{\rm max} (1) \eqd
\max \{ UM^{\{1\}}, (1-U)M^{\{2\}} \} ,\]
where $U$, $M^{\{1\}}$ and $M^{\{2\}}$ are independent
random variables, $U$ is uniform on $(0,1)$, 
and $M^{\{1\}}$ and $M^{\{2\}}$ have the max-Dickman distribution
as given by (\ref{maxdick}).
\end{theorem}

 We will derive Theorem \ref{dickman} from a new result
 on the limiting distribution of the length of the longest edge in the ONG
 on uniform random points in $(0,1)^d$, which is of some independent interest: see Theorem \ref{ongmax} below.

As promised, we now give a characterization of the limits $Q(1,\alpha)$, $\alpha \geq 1$,
arising in the $d=2$ case of Theorem \ref{mainth}.
First we define random variables $\tJal$, $\alpha >1/2$, with
$\Exp[\tJal]=0$ and $\Exp[ \tJal^2]<\infty$.
Define 
   $\tJo$ by the fixed-point equation
\begin{align}
\label{0919x}
\tJo \eqd \min \{ U, 1-U \} + U \tJo^{\{1\}} + (1-U) \tJo^{\{2\}}
+ \frac{1}{2} U \log U + \frac{1}{2} (1-U) \log (1-U) ,
\end{align}
and for $\alpha \in (1/2,\infty) \setminus \{1\}$, define $\tJal$ by the fixed-point equation
\begin{align}
\label{0919z}
\tJal \eqd \min \{ U, (1-U) \}^\alpha + U^\alpha \tJal^{\{1\}}
+ (1-U)^\alpha \tJal^{\{2\}} +  \frac{2^{-\alpha}}{\alpha-1}
\left( U^\alpha +(1-U)^\alpha -1 \right).
\end{align}
In each of these two equations (and subsequently), $Y^{\{1\}}$ and $Y^{\{2\}}$ denote
independent copies of the random variable $Y$, and $U$ denotes a uniform
random variable on $(0,1)$ independent of the other random variables
on the right-hand side of the equation.

Note that (\ref{0919x}) and (\ref{0919z})
define unique square-integrable mean-zero solutions (see e.g.~Theorem 3 of R\"osler \cite{rosler0}),
and hence the distributions of $\tJo$ and $\tJal$ are uniquely defined. 
 Moments of $\tJal$
  can be calculated recursively from (\ref{0919x}) and (\ref{0919z}); see
\cite{ong} for some information on the first few
moments of $\tJo$, for example. From these  moments
one can deduce that $\tJal$, $\alpha >1/2$ is not Gaussian.

Now we can define random variables $\tHal$, $\tGal$, again
with zero mean and finite variance.
Define $\tHo$ by
\begin{align*}
\tHo \eqd U \tJo + (1-U) \tHo + \frac{U}{2} +
\frac{1}{2} U \log U  
+ \frac{1}{2} (1-U) \log (1-U). 
\end{align*} 
For $\alpha \in (1/2,\infty) \setminus \{1\}$, define $\tHal$ by
\begin{align*}
\tHal \eqd U^\alpha \tJal
+ (1-U)^\alpha \tHal + U^\alpha \left( 1+ \frac{2^{-\alpha}}{\alpha -1} \right)
+ ((1-U)^\alpha-1) \left( \frac{1}{\alpha} + \frac{2^{-\alpha}}{\alpha (\alpha-1)} \right).
\end{align*}
Define $\tGo$ by 
\begin{align}
\label{0926a}
\tGo \eqd U \tHo^{\{1\}} 
+(1-U) \tHo^{\{2\}} +\frac{1}{4} + \frac{1}{2} U \log U +\frac{1}{2} (1-U) \log (1-U).
\end{align}
Finally, for $\alpha \in (1/2,\infty) \setminus \{1\}$, define 
$\tGal$ 
by
\begin{align}
\label{0926p}
\tGal   \eqd {} &
  U^\alpha \tHal^{\{1\}} 
+(1-U)^\alpha \tHal^{\{2\}}+ \left( U^\alpha +(1-U)^\alpha \right) \left(
\frac{1}{\alpha} + \frac{2^{-\alpha}}{\alpha (\alpha -1)} \right) \nonumber\\ 
& 
- \frac{2}{\alpha (\alpha +1)} \left( 1 +\frac{2^{-\alpha}}{\alpha -1 } \right).
\end{align} 
Once again,
the distributions of   $\tHal$  and $\tGal$ are uniquely defined. It is the distribution
of $\tGal$ $(\alpha \geq 1)$ as defined by
(\ref{0926a}) or (\ref{0926p}) that appears
in the $d=2$ case  of Theorem \ref{mainth}.

In the remainder of this paper, we prove Theorems \ref{llnthm}, \ref{mainth} and \ref{dickman}.
First, in Section \ref{secong} we discuss
the ONG,
 which we use to deal with the boundary
effects in the MDST, and prove some new results, which are of some independent interest. 
In 
Section \ref{middle},
we apply general results
of Penrose \cite{mpbern,mpejp} (see also \cite{mdp05})
to prove a law of large numbers and central limit theorem for
the total weight of the MDST away from the boundary.  
In Section \ref{bdry}
we deal with the boundary effects themselves. Then in Section \ref{secdick}
we prove Theorem \ref{dickman}. Finally, we complete the proofs of
Theorem \ref{mainth} in
Section \ref{totallength} and Theorem \ref{llnthm} in Section \ref{prflln}.

Throughout the sequel
we make repeated use of {\em Slutsky's theorem} (see, e.g., Durrett \cite{durrett}, p.~72), which says
that for sequences of random variables $(X_n)$, $(Y_n)$ such that $X_n \tod X$ and
$Y_n \toP 0$ as $n\to\infty$, we have $X_n+Y_n \tod X$ as $n \to \infty$. (Here and subsequently
`$\toP$' denotes convergence in probability.)

\section{The on-line nearest-neighbour graph}
\label{secong}

In this section we describe the on-line nearest-neighbour graph
that we use 
to analyse the boundary effects in the total
weight of the MDST under $\pos$. Some of the
results that
we will require are present in \cite{ong} and \cite{ong2}, but we will also need some
new results, which we prove in this section.

Let $(\bY_1,\bY_2,\ldots)$ be a sequence of vectors
in $\R^d$, and for $m \in \N$ set
$\Y_m := (\bY_1,\ldots,\bY_m)$.
The on-line nearest-neighbour graph (ONG)
on sequence
 $\Y_m$ is constructed by joining each point 
 after the first of $\Y_m$ by a
 directed edge
 to its (Euclidean)
 nearest neighbour
amongst those points that precede it in the sequence. That is, for $i=2,\ldots,m$
we include the edge $(\bY_i,\bY_j)$ where $j \in \{1,\ldots,i-1\}$
is such that
 \[ 
\| \bY_j -\bY_i\|_d
= \min_{1 \leq k < i} \| \bY_k - \bY_i \|_d, \]
arbitrarily breaking any ties.
  
In this way we obtain the ONG on $\Y_m$,
denoted $\mathrm{ONG} ( \Y_m )$ and which, ignoring
directedness of edges,
 is a tree rooted at $\bY_1$. 
 Denote the total power-weighted
 edge-length with exponent $\alpha>0$ of 
 $\mathrm{ONG} ( \Y_m )$ by
$\O^{d,\alpha} ( \Y_m )$, that is
\[ \O^{d,\alpha} (\Y_m) := \sum_{i=2}^m \min_{1 \leq j < i} \| \bY_i - \bY_j \|_d^\alpha ;\]
 when $\Y_m$ is random, we denote the centred  version by
$\tO^{d,\alpha} ( \Y_m ) := \O^{d,\alpha} ( \Y_m
) - \Exp [ \O^{d,\alpha} ( \Y_m )
]$.

 Our primary interest is the case where
 $\Y_m$ is a sequence of uniform random vectors on the $d$-cube. 
 Let $d \in \N$.
Let $(\bU_1,\bU_2,\ldots)$ be a sequence of independent uniform random vectors
in $(0,1)^d$. For $m \in \N$, set
$\U_m := (\bU_1, \ldots, \bU_m)$. We then consider $\mathrm{ONG} (\U_m)$.

We also consider the ONG defined on a Poisson number of points.
Let
$(N(t))_{t\geq 0}$ be the counting process of a homogeneous Poisson process
of unit rate in $(0,\infty)$, independent of $(\bU_1,\bU_2,\ldots)$.
Thus $N(n)$ is a Poisson random
variable with mean $n$.
With $\U_m$ as defined above set 
$\Pi_n = \U_{N(n)}$;
we then consider $\mathrm{ONG} (\Pi_n)$
Note that the points
of the sequence $\Pi_n$ then constitute a homogeneous  
Poisson point process
of intensity $n$ on $(0,1)^d$.

 We need the following result, which is contained in Theorem 2.1 of \cite{ong2}.

\begin{lemma}
\label{ongl}
Suppose $d \in \N$.
\begin{itemize}
\item[(i)] For $\alpha \in (0,d/2)$,
there exists a constant $C \in (0,\infty)$ such that for all $n \geq 1$
\[ \Var [\tO^{d,\alpha} (\Pi_n)] \leq C n^{1-(2\alpha/d)}.\]
\item[(ii)] For $\alpha=d/2$, 
 there exists a constant $C \in (0,\infty)$ such that for all $n \geq 1$
\[ \Var [ \tO^{d,d/2} (\Pi_n)] \leq C \log (1+n).\]
\end{itemize}
\end{lemma}

The following result is contained in Theorem 2.2 of \cite{ong2}, with
Theorem 2.2 of \cite{ong} used to deduce the final statement about the $d=1$ case.

\begin{lemma}
\label{onngthm}
Suppose $d \in \N$ and $\alpha>d/2$. Then 
there exists a mean-zero random variable $Q(d,\alpha)$ such that as $n\to \infty$
\begin{align*}
 \tO^{d,\alpha} (\Pi_n) \tod Q(d,\alpha).\end{align*}
In particular, $Q(1,\alpha) = \tGal$ for $\alpha \geq 1$, where
$\tGal$ has distribution given by (\ref{0926a}) or (\ref{0926p}).
\end{lemma}

In order to deduce
 Theorem \ref{dickman}, we use the following result on the length of longest edge of the ONG
 on uniform random points in $(0,1)^d$, which adds to the
 analysis of the  ONG given in \cite{boll,mdp,ong,llns,ong2}.
For a sequence $\Y_m =(\bY_1,\ldots,\bY_m)$ of points in $\R^d$, let $\Om^d (\Y_m)$ denote the length of the longest
 edge in the ONG on $\Y_m$:
 \[  \Om^d (\Y_m) := \max_{2 \leq i \leq m} \min_{1 \leq j < i} \| \bY_i - \bY_j \|_d .\] 
 For $d=1$, 
where $\U_n = (U_1,\ldots, U_n)$ and
 $\Pi_n = (U_1,\ldots,U_{N(n)})$ 
 for $U_1,U_2,\ldots$
 independent
 uniform random variables on $(0,1)$, we set
 $\U_n^0 := (0,U_1,\ldots,U_n)$, i.e.~$\U_n^0$ is $\U_n$ but with
 an initial point placed at the origin, and similarly
 $\Pi_n^0 := (0,U_1,\ldots,U_{N(n)} )$.
  
\begin{theorem}
\label{ongmax}
Let $d \in \N$.
\begin{itemize}
\item[(i)] There exists a random variable $Q_{\rm max} (d)$ such that
as $n \to \infty$
\[ \Om^d ( \U_n ) \tod Q_{\rm max} (d); ~~~ \Om^d (\Pi_n ) \tod Q_{\rm max} (d).\]
\item[(ii)] When $d=1$, we have in particular that
\begin{equation}
\label{pp0}
 Q_{\rm max} (1) \eqd \max \{ U M^{\{1\}}, (1-U) M^{\{2\}}\} ,\end{equation}
where $U$, $M^{\{1\}}$, $M^{\{2\}}$ are independent, $U$ is uniform
on $(0,1)$ and $M^{\{1\}}$, $M^{\{2\}}$ are max-Dickman random variables 
as given by (\ref{maxdick}). Also as $n \to \infty$
\begin{align}
\label{mm1}
 \Om^1 (\U^0_n) \tod M; ~~~ \Om^1 (\Pi^0_n)  \tod M ,\end{align}
where $M$ is a max-Dickman random variable 
as given by (\ref{maxdick}).
\end{itemize}
\end{theorem}
\proof 
First we prove part (i). With probability $1$,
for all $n$, $0 \leq \Om^d (\U_n) \leq d^{1/2}$   and
$\Om^d (\U_{n+1}) \geq \Om^d (\U_n)$. Hence $\Om^d ( \U_n ) \to Q_{\rm max} (d)$
a.s., as $n \to \infty$,
for some $Q_{\rm max} (d)$. Then by the coupling of $\Pi_n$ and $\U_n$
and the fact that $N(n) \to \infty$ a.s., we have that with this coupling
 $\Om^d (\Pi_n)$ converges to the same $Q_{\rm max}$ a.s.~and hence in 
 distribution (regardless of the coupling), 
completing the proof of part (i).

We now prove part (ii) of the theorem, and so take $d=1$. First we
prove (\ref{mm1}).
Again by the coupling of $\Pi_n$ and $\U_n$, it suffices 
to prove that $\Om^1(\U^0_n) \to M$ a.s. as $n \to \infty$.
The following argument is related to the proof of Theorem 2 of \cite{rooted}.

 An {\em upper record value} 
in the sequence $X_1, X_2,X_3, \ldots$ is a
value $X_i$ which exceeds $\max\{X_1$, $\ldots$, $X_{i-1} \}$
(the first value $X_1$ is also included as a record value). Let $j(1), j(2),\ldots$ be the values of
$i \in \N$ such that $U_i$ is an upper record in the sequence
$(U_1,U_2,\ldots)$, arranged in increasing order so that $1=j(1) < j(2) < \cdots$.
Let $R_n := \max\{ k: j(k) \leq n \}$
 be the number of record values in the sequence
$\U_n = (U_1,\ldots,U_n)$.

A record $U_i$ has by definition
 no preceding point
in the sequence $\U_n$ to its right in the unit
interval, and hence (in the ONG on $\U^0_n$)
 must be joined
to its nearest neighbour to the left amongst those points
already present, which is necessarily the previous record value when $i>1$, or $0$
in the case of $U_1$. 
Then each non-record $U_i$ lies in an interval between
a record value and its nearest neighbour to the left, and hence
gives rise to a shorter edge than that from some record value. Thus
\begin{align}
\Om^1 (\U_n^0  ) = \max_{1 \leq i \leq R_n} \left\{ U_{j(i)} - U_{j(i-1)} \right\},
\label{1108}
\end{align}
where we set $j(0) :=0$ and $U_0 :=0$.
For $i \in \N$ set
$$
V_i := \frac{ 1 - U_{j(i)} }{1-U_{j(i-1)}}.
$$
It is not hard to see that $V_1, V_2, \ldots $ are
mutually independent and each is
 uniformly distributed over $(0,1)$.
Therefore, setting
\[ M :=   \max\{ 1- V_1, V_1(1- V_2), V_1V_2(1 - V_3), V_1V_2V_3(1-V_4),\ldots \}  ,\]
we obtain
\begin{align}
\label{Mdef}
M & =   \max \{ 1-V_1 , V_1 \max \{ 1-V_2, V_2 (1-V_3), V_2 V_3(1-V_4),\ldots \}\} \nonumber\\
& =   \max \{ 1-V_1, V_1 M' \}, 
\end{align}
where $M':=\max \{ 1-V_2, V_2 (1-V_3), V_2 V_3(1-V_4),\ldots \}$
has the same distribution as $M$ and is independent of $V_1$.
Hence $M$   has the max-Dickman
distribution as given by (\ref{maxdick}).
Furthermore, with the convention that an empty product is $1$,
\begin{equation}
(1 - V_i)
\prod_{k=1}^{i-1} V_k 
 =
 \frac{U_{j(i)} - U_{j(i-1)} }{1 - U_{j(i-1)}} 
\prod_{k=1}^{i-1} \left( \frac{1-U_{j(k)}}{1-U_{j(k-1)}} \right)
= U_{j(i)} - U_{j(i-1)} ,
\label{1108a}
\end{equation}
for $k \in \N$.
Also,
$R_n \to \infty$ almost surely as $n \to
\infty$. Hence by (\ref{1108}), (\ref{Mdef}) and
(\ref{1108a}), 
\[ \Om^1 (\U^0_n) = \max_{1 \leq i \leq R_n} \left\{
(1-V_i) \prod_{k=1}^{i-1} V_k \right\}
\to \max_{i \geq 1} \left\{
(1-V_i) \prod_{k=1}^{i-1} V_k \right\} = M ,\]
where the convergence is almost sure.
This proves (\ref{mm1}).

To complete the proof of part (ii) of the theorem, we 
need to prove (\ref{pp0}). Conditioning on $U = U_1$ 
and the number of points of $(U_2,U_3, \ldots, U_n)$
that fall in each of the two intervals $(0,U)$, $(U,1)$,
we obtain by scaling that
\begin{equation}
\label{pwp1}
 \Om^1 (\U_n) \eqd \max \{ U  \Om^1 ( \U^0_{L} )
, (1-U) \Om^1 ( \tilde \U^0_{n-1-L}) \} ,\end{equation}
where in
 the right-hand expression $\tilde \U^0_m = (0, \tilde U_1, \tilde U_2, \ldots, \tilde U_m)$,
$L \sim {\rm Bin} (n-1, U)$, and $U, U_1, U_2, \ldots, \tilde U_1, \tilde U_2, \ldots$
are independent uniform random variables on $(0,1)$. Here $L$ and $n-1-L$ both
tend to infinity
a.s.~as $n \to \infty$,
and
$\Om^1 ( \U^0_{L} )$ and $\Om^1 ( \tilde \U^0_{n-1-L})$
are independent given $L$. Thus by
 (\ref{mm1}) we have that
$\Om^1 ( \U^0_{L} )$ and $\Om^1 ( \tilde \U^0_{n-1-L})$
converge in distribution to independent
copies of the max-Dickman variable $M$.
Then (\ref{pwp1}) and
the fact that $Q_{\rm max} (1)$ is the distributional limit
of $\Om^1 (\U_n)$ 
yields (\ref{pp0}). $\square$

\section{Limit theorems away from the boundary}
\label{middle}

In this section we prove a law of large numbers and central
limit theorem for the total power-weighted length of the MDST edges
from points that are not too close to the base of the unit $d$-cube.
To do this, we employ some general results of Penrose \cite{mdp05,mpbern,mpejp}.

Recently, notions of {\em stabilizing} functionals
of point sets have proved to be a useful basis for
a general methodology for establishing  limit
theorems
for functionals 
of random
point sets in $\R^d$. See for example \cite{mdp,mpbern,mpejp,penyuk1,penyuk2}.
To prove the law of large numbers (Lemma \ref{llnlem})
and central limit theorem (Lemma \ref{CLT})
in this section, we make use of the general
results on convergence of random measures in geometric
probability given in \cite{mdp05,mpbern,mpejp}. These two lemmas will then
form two of the ingredients for two of our main results,
Theorems \ref{llnthm} and \ref{mainth}.

We use the following notation. Let $d \in \N$.
Let $\X \subset \R^d$ be finite. For 
constant $a>0$, and ${\bf y} \in
\R^d$, let ${\bf y}+a\X$ denote the transformed set $\{ {\bf y} + a\bx
: \bx \in \X\}$. 
For $\bx \in \R^d$ and $r>0$, let $B(\bx;r)$ be the closed
 Euclidean $d$-ball
with centre $\bx$ and radius $r$.
For 
bounded measurable $R \subset \R^d$ let $|R|$ denote the 
 $d$-dimensional
Lebesgue measure 
of $R$. Write $\0$ for
 the origin of $\R^d$.

For $\alpha>0$, define the $[0,\infty)$-valued
function on finite non-empty $\X \subset \R^d$ and $\bx \in \X$:
\begin{align}
\label{0606vv}
\xi(\bx;\X) := d_*(\bx;\X)^\alpha,\end{align}
and set $\xi(\bx; \emptyset):=0$ for any $\bx$.
Then $\xi$ is translation invariant (that is $\xi( \by + \bx ; \by + \X) = \xi(\bx ; \X)$ for all
$\by \in \R^d$, all finite $\X \subset \R^d$ and $\bx \in \X$)
and homogeneous of order $\alpha$ (that is for any $r>0$, $\xi ( r \bx ; r \X) = r^\alpha \xi( \bx;\X)$
for all finite $\X \subset \R^d$ and $\bx \in \X$). 
For $\X \subset \R^d$ and $\bx \in \R^d$, write $\X^\bx$ for $\X \cup \{ \bx\}$.
If $\bx \notin \X$, we abbreviate notation
to $\xi(\bx;\X)=\xi(\bx ; \X^\bx)$. The above definitions
extend naturally to infinite but locally
finite sets $\X$ (as in \cite{mpbern}).

Let 
\begin{align}
 \LL^{d,\alpha} (\X ; R) : = \sum_{\bx \in \X \cap R} \xi (\bx ; \X) 
\label{0714a}
\end{align}
be the translation invariant functional 
defined on all 
finite point sets $\X \subset \R^d$ and all Borel sets $R \subseteq \R^d$ induced by the function $\xi$. 
Then $\LL^{d,\alpha}(\X;R)$ is the total power-weighted
length of the edges of the MDST on $\X$ originating
from points in the region $R$. It is this functional that interests us here. When $\X$ is random,
 set $\tilde
\LL^{d,\alpha} (\X;R):= \LL^{d,\alpha} (\X;R) - \Exp[  \LL^{d,\alpha} (\X;R)]$. Note that with our previous 
notation, $\LL^{d,\alpha} (\X)=\LL^{d,\alpha} (\X; (0,1)^d)$ for $\X \subset (0,1)^d$.

Fix $\eps \in (0,1/d)$ (small). Let $(g_n)_{n > 0}$ be
such that $g_n \in (0,1)$ and
$g_n = \Theta (n^{\eps-(1/d)})$ as $n \to \infty$, where by
$a(n) = \Theta (b(n))$ as $n \to \infty$
we mean \[
0 < \liminf_{n \to \infty} \frac{a(n)}{b(n)} \leq
 \limsup_{n \to \infty} \frac{a(n)}{b(n)} < \infty .\]
Given $g_n$,
we introduce the 
family $(\Gamma_n)_{n \geq 1}$
of Borel subsets of $\R^d$
given by 
\begin{align}
\label{77a}
\Gamma_n := (0,1)^{d-1} \times (g_n,1),
\end{align} i.e.~$\Gamma_n$ is the unit $d$-cube without a thin strip
at the base (in the $\be_d$ sense). Note that the limiting set $\cup_{n \geq 1} \Gamma_n = (0,1)^d$.
Later on, in Section \ref{totallength},
we will make a more specific choice for $g_n$.
For $n \geq 1$, locally finite $\X \subset \R^d$ and $\bx \in \X$ we define the scaled-up
 version of $\xi$ 
restricted to $\Gamma_n$ by
\[ \xi_n (\bx ;\X) := \xi (n^{1/d} \bx ; n^{1/d} \X ) \1_{\Gamma_n} (\bx).\]
Then, from (\ref{0714a}) 
\begin{align}
\label{966}
 \LL^{d,\alpha} (\X ; \Gamma_n ) = \sum_{\bx \in \X}
\xi (\bx ; \X) \1_{\Gamma_n} (\bx)
= n^{-\alpha/d} \sum_{\bx \in \X} \xi_n ( \bx ; \X),\end{align}
using the fact that $\xi$ as given by (\ref{0606vv})
is homogeneous of order $\alpha$.
 We employ the following notion of stabilization (see \cite{mpbern,mpejp}).
 
\begin{definition}
\label{stabdef} 
For any locally finite $\X \subset \R^d$ and Borel region $A \subseteq \R^d$,
 define $R_\xi(\0;\X,A)$ (called the {\em
radius of stabilization} for $\xi$ at $\0$ with respect to $\X$ and $A$) to be the smallest
integer $r\geq 0$ such that
\[ \xi(\0; (\X \cap B(\0;r))\cup \Y) = \xi (\0; \X \cap B(\0;r) ),\]
for all finite $\Y \subset A \setminus B(\0;r)$. If no such $r$ exists,
set $R_\xi (\0;\X,A)=\infty$. 
\end{definition}
When $A$ is all of $\R^d$, we write $R_\xi (\0;\X)$ for $R_\xi (\0;\X,\R^d)$.

\subsection{Law of large numbers}
\label{llnproof}

We will apply a Poisson point process analogue
of the law of large numbers Theorem 2.1 of \cite{mpbern}.
As mentioned on p.~1130 of \cite{mpbern}, such a Poisson-sample
result follows by similar arguments to the proofs
in \cite{mpbern}; in fact such a result is stated
and proved as Theorem 2.1 in \cite{mdp05}. It is this latter result that we will use in
this section.

Let $\H_1$ denote 
a homogeneous Poisson point process of unit intensity on $\R^d$.
 Our law of large numbers
result for this section is the following.

\begin{lemma}
\label{llnlem}
Suppose $d \in \{2,3,\ldots\}$ and $\alpha >0$.
As $n \to \infty$ we have 
\begin{align}
\label{77}
n^{(\alpha/d)-1} \LL^{d,\alpha} (\Po_n ; \Gamma_n ) 
 \to \Exp [ \xi (\0 ; \H_1 )]
 =
2^{\alpha/d} v_d^{-\alpha/d} \Gamma( 1+(\alpha/d) )
, \end{align}
where the convergence is in $L^2$, and
 $v_d$ is given by (\ref{0818c}).
 \end{lemma}

The statement (\ref{77}) will follow
 from Theorem 2.1 of \cite{mdp05} applied to our
functional $\xi$ as defined at (\ref{0606vv}), using (\ref{966}). Thus we need to verify the conditions of
Theorem 2.1 of \cite{mdp05}: (a) that $R_\xi(\0;\H_1)$ is almost surely finite; and
(b) that there exists
some $p>2$ such that the following two moments conditions hold:
\begin{align}
\label{mom1}
\sup_{n \geq 1; ~\bx \in (0,1)^d} \Exp [ \xi_n (\bx ; \Po_n )^p ] < \infty,~~\textrm{and}\\
\sup_{n \geq 1; ~\bx,\by \in (0,1)^d} \Exp [ \xi_n (\bx ; \Po^\by_n )^p ] < \infty.
\label{mom2}
\end{align}
The next two lemmas take care of this.
 
\begin{lemma}
\label{stab2}
For $\xi$ given by (\ref{0606vv}), 
the radius of stabilization $R_\xi(\0;\H_1)$ as defined in Definition \ref{stabdef}
is almost surely finite.
\end{lemma}
\proof
Let $R=d_*(\0;\H_1)$. Then $R$ is finite almost surely. For
 any $\ell>R$ we have that $\xi(\0;(\H_1 \cap B(\0;\ell))\cup {\cal Y})=R^\alpha$,
for any finite ${\cal Y} \subset \R^d \setminus B(\0;\ell)$. Thus taking $R_\xi(\0;\H_1)$
to be the smallest integer greater than $R$, $R_\xi(\0;\H_1)$ is 
almost surely finite. $\square$

\begin{lemma}
\label{momlem}
Suppose $d \in \{2,3,\ldots\}$ and $\alpha>0$.
Then for $(\Gamma_n)_{n \geq 1}$
as given at (\ref{77a}) and $\xi$ as given by (\ref{0606vv})
 the moments conditions (\ref{mom1}) and (\ref{mom2}) hold
 for any $p >0$.
\end{lemma}
\proof
We have from the definition of $\xi_n$ and (\ref{0606vv}) that
\begin{align}
 \sup_{\bx \in (0,1)^d} \Exp [ \xi_n (\bx ; \Po_n)^p] = \sup_{\bx \in \Gamma_n}
\Exp [ \xi (n^{1/d}\bx ; n^{1/d} \Po_n)^p ]
= \sup_{\bx \in \Gamma_n} \Exp[ d_* (n^{1/d} \bx ; n^{1/d} 
\Po_n )^{\alpha p} ].
\label{0118a}
\end{align}
For $d \in \{2,3,\ldots\}$,
$\bx \in \Gamma_n$ and $r>0$, define the region in the scaled-up
space $(0,n^{1/d})^d$
\begin{align}
\label{0612a}
A^d_n(\bx,r):= B(n^{1/d} \bx;r) \cap (0,n^{1/d})^d \cap \{ \by \in \R^d : \by \pos n^{1/d} \bx \}.
\end{align}
For $\bx \in \Gamma_n$,
define the  variables 
$\zeta_n^{(1)} (\bx) := d_*(n^{1/d}\bx;n^{1/d}\Po_n)
{\bf 1}_{\{d_*(\bx;\Po_n) \leq g_n\}}$ and
$\zeta_n^{(2)} (\bx) := d_*(n^{1/d}\bx;n^{1/d}\Po_n)
{\bf 1}_{\{d_*(\bx;\Po_n) > g_n\}}$.
For $t \geq 0$, 
\[
\Pr ( 
\zeta_n^{(1)} (\bx) > t ) = 
\Pr ( \{ 
d_*(n^{1/d}\bx;n^{1/d}\Po_n) 
>t \}
\cap \{ d_*(n^{1/d}\bx;n^{1/d}\Po_n) \leq 
n^{1/d}g_n\} ).
\]
This probability is clearly zero
unless $t < n^{1/d} g_n$, in which case,
 by the definition of $\Gamma_n$ the  region $A_n^d(\bx,t)$ does not
touch the  hyperplane $\{x_d = 0\}$, so that $|A_n^d(\bx,t)|
 \geq 2^{-d} v_d t^d$, 
 where
$v_d$ is the volume of the unit $d$-ball given by (\ref{0818c}). 
Hence for all $t \geq 0$,
\[
\Pr ( 
\zeta_n^{(1)} (\bx) 
> t ) \leq 
\exp(  - 2^{-d} v_d t^d ) 
\]
so that for all $n$ and all $\bx \in \Gamma_n$, 
$\zeta_n^{(1)} (\bx) $
is stochastically dominated by a variable with cumulative
distribution function $F(t) = 1 - 
\exp(  - 2^{-d} v_d t^d )$, $t \geq 0$. Such a variable has finite $(\alpha p)$-th moment.

Also, for all $n$ and all $\bx \in \Gamma_n$, 
the random variable 
$ \zeta_n^{(2)} (\bx)$
 is bounded by
 the random variable 
$d^{1/2}n^{1/d}
\1_{ \{d_*(n^{1/d}\bx;n^{1/d}\Po_n)> n^{1/d} g_n \}}$, so that  
\begin{eqnarray*}
\Exp [
(\zeta_n^{(2)} (\bx) 
)^{\alpha p} ]
\leq d^{\alpha p/2}n^{\alpha p/d}
\Pr ( d_*(n^{1/d}\bx;n^{1/d}\Po_n)> n^{1/d} g_n ) 
\\
\leq d^{\alpha p/2}n^{\alpha p/d}
\exp( - |A_n^d(\bx,n^{1/d} g_n)|)
\leq d^{\alpha p/2}n^{\alpha p/d}
\exp( - 2^{-d} v_d (n^{1/d} g_n)^d)
\end{eqnarray*}
and since $n^{1/d} g_n = \Theta(n^{\varepsilon})$, this upper bound
is bounded in $n$. 
Thus the $(\alpha p)$-th moment of 
$\zeta_n^{(2)} (\bx)$ is bounded uniformly over
 all $n$ and all $\bx \in \Gamma_n$. 
Combined with the earlier uniform moment bound for
$\zeta_n^{(1)} (\bx)$ and (\ref{0118a}), this yields
(\ref{mom1}).

For (\ref{mom2}), note that for any
 $\bx \in \Gamma_n, \by \in (0,1)^d$
\[ \xi_n ( \bx; \Po_n^\by ) = d_* ( n^{1/d} \bx ; n^{1/d} ( \Po_n \cup \{ \by \} ) )^\alpha
\leq d_* ( n^{1/d}\bx ; n^{1/d}\Po_n)^\alpha +
 \1_{\{ \Po_n \subset \Gamma_n \}} n^{\alpha/d} d^{\alpha/2} .\]
Moreover, $
 \xi_n ( \bx; \Po_n^\by ) $
is zero for $\bx \in (0,1)^d \setminus \Gamma_n$.
Thus 
\[ \sup_{\bx, \by \in (0,1)^d} \Exp \left[
\xi_n (\bx ; \Po_n^\by) ^p \right] 
\leq \sup_{\bx \in (0,1)^d} \Exp \left[ \xi_n (\bx ; \Po_n)^p \right] +
 \Pr ( \Po_n \subset \Gamma_n ) n^{\alpha p/d} d^{\alpha p/2} ,\]
so that (\ref{mom1}) implies (\ref{mom2})
since $\Pr (\Po_n \subset \Gamma_n) = \exp(-n g_n)  $. $\square$ \\

\noindent
{\bf Proof of Lemma \ref{llnlem}.} 
From Theorem 2.1 of \cite{mdp05}, with (\ref{966})
 and Lemmas \ref{stab2} and
\ref{momlem}, we obtain the convergence statement in (\ref{77}).
It remains to prove the final equality (\ref{77}). We have, for $s \geq 0$
\begin{align*}
\Pr\left( \xi (\0;
\H_1 ) >s \right)   =  
\Pr \left( \H_1 \cap
\{ \bx \in \R^d : \bx \pos \0 \} \cap B(\0;s^{1/\alpha}) = \emptyset \right)  = 
\exp{(- (v_d/2) s^{d/\alpha} )}.
\end{align*}
Hence,
 \[
  \Exp \left[ \xi (\0; \H_1 ) \right]  =
\int_0^{\infty} \Pr \left( \xi (\0; \H_1 )
> s \right) \ud s
= \int_0^\infty \exp ( - (v_d/2) s^{d/\alpha} ) \ud s,\]
which by the change of variables $y = (v_d/2)s^{d/\alpha}$ is the same as
\[ \frac{\alpha}{d} 2^{\alpha/d} v_d^{-\alpha/d} \int_0^\infty y^{(\alpha/d)-1} \exp (-y) \ud y
= \frac{\alpha}{d} 2^{\alpha/d} v_d^{-\alpha/d} \Gamma ( \alpha /d),\]
by Euler's Gamma integral (see e.g.~6.1.1 in \cite{as}). The desired equality
 now follows from the functional relation
$x \Gamma (x) = \Gamma (1+x)$ (see 6.1.15 in \cite{as}). $\square$

\subsection{Central limit theorem} 
\label{ltot}

We again consider $\LL^{d,\alpha}(\Po_n ; \Gamma_n)$ as given by (\ref{966}).
In this section we aim to prove a central limit theorem complementing
 the law of large numbers of Section \ref{llnproof}. 
 This time, we will apply Theorems 2.1 and 2.2 of
\cite{mpejp} to give the following result.

\begin{lemma} \label{CLT} 
Let $d \in \{2,3,\ldots\}$ and $\alpha >0$.
There exists a constant $s_\alpha \in [0,\infty)$, 
not depending on the choice of $\eps$ or the sequence $g_n$,
such that, as $n \to \infty$,
\[ \lim_{n \to \infty} \left( n^{(2\alpha/d)-1} \Var \left[
\LL^{d,\alpha} ( \Po_n ; \Gamma_n ) \right] \right)
= \lim_{n \to \infty} \left( n^{-1} \Var \left[ \sum_{\bx \in \Po_n}
\xi_n (\bx ; \Po_n) \right] \right) = s_\alpha^2,\]
and
\[
n^{(\alpha/d)-(1/2)} \tilde \LL^{d,\alpha}
 ( \Po_n ; \Gamma_n ) 
 \tod \NN ( 0, s_\alpha^2 ).\]
\end{lemma}
\proof
In order to prove this lemma, we need to verify the
conditions of Theorems 2.1 and 2.2 
 of \cite{mpejp} 
(see also Theorem 2.2 and 2.3 of \cite{mdp05}) for
our function $\xi$ as given by (\ref{0606vv}). In addition
to the moments conditions (\ref{mom1}), (\ref{mom2}) (as shown to
hold in Lemma \ref{momlem}), we need to demonstrate the following
additional stabilization conditions:
\begin{align}
\label{2a}
\Pr ( R_\xi(\0; \H_1^\bz ) < \infty ) =1,\end{align}
for all $\bz \in \R^d$; and
\begin{align}
\label{2b}
\limsup_{s \to \infty}
s^{-1} \log \left( \sup_{n \geq 1; ~ \bx \in \Gamma_n} \Pr ( R_\xi( n^{1/d} \bx ; n^{1/d} \Po_n 
, n^{1/d} (0,1)^d ) > s) \right)
< 0 .\end{align}
Condition (\ref{2a}) requires that the radius of stabilization
is almost surely finite on the addition of an arbitrary extra
point to $\H_1$, and condition (\ref{2b})
requires exponential decay of the tail of the radius of stabilization.

Given Lemma \ref{stab2}, (\ref{2a}) is clear, since
with probability $1$ the addition of any extra point
$\bz \in \R^d$ to $\H_1$ can only decrease the radius of stabilization at $\0$.

We need to prove (\ref{2b}).
Let
  $|A^d_n( \bx , r )| $ be defined by (\ref{0612a}), and  
for $\bz=(z_1,z_2,\ldots,z_d) \in n^{1/d} \Gamma_n$, let $m(\bz):= z_d$,
 the last component
of $\bz$. 
For $d \geq 2$, we claim that
there are finite constants $C_d > 0$ and $n_0 \geq 1$ such that
\begin{align}
\label{0612b}
  |A^d_n( \bx , r )| \geq 
C_d r^{d-1}  ~
& \textrm{if} ~ r \in (1,d^{1/2}n^{1/d}], 
\end{align} 
for all $n \in \N$ with $n \geq n_0$, and any $\bx \in \Gamma_n$. 

We verify  the claim (\ref{0612b}). Take $n_0$ such that
for all $n \geq n_0$ we have $n^{1/d}g_n \geq 1$.
Then for $n \geq n_0$,
suppose $r \in (1,d^{1/2} n^{1/d} ]$.
For a lower bound on the volume of $A^d_n( \bx , r )$,
consider $\bx = (0,0,\ldots,0,m(\bx))$, the `worst case'.
 Let $h_\bx$ denote
the hyperplane $\{ \by \in n^{1/d} \Gamma_n : m(\by)=m(n^{1/d} \bx)\}$.
Let $r' := d^{-1/2}r$, so $r' \leq n^{1/d}$.
 Then let $\bw_1,\bw_2,\ldots,\bw_{d-1}$
denote the $d-1$ points of $h_\bx$
$(r',0,0,\ldots,0,m(n^{1/d} \bx))$, $(0,r',0,\ldots,0,m(n^{1/d} \bx))$,
\ldots
$(0,0,\ldots ,0,r',m(n^{1/d} \bx))$, and let 
$\bw_0$ denote the point $(0,0,\ldots,0,m(n^{1/d} \bx ) -1)$.
 Then 
since $\bx \in \Gamma_n$,
the $d$-dimensional `right pyramid' defined
by vertices $\bw_0,n^{1/d} \bx,\bw_1,\ldots,\bw_{d-1}$ is
 contained within both $(0,n^{1/d})^d$ and
the half-ball
$B(n^{1/d} \bx;r) \cap \{ \by \in \R^d : \by \pos n^{1/d} \bx \}$. 
The volume of this `pyramid'
is
$ d!^{-1}  (r')^{d-1}$. 
This gives
a lower bound
for $|A^d_n( \bx , r )|$,  and 
 (\ref{0612b}) holds as claimed.

To prove (\ref{2b}), note that $n^{1/d} \Po_n$ is a homogeneous Poisson point process
of unit intensity on $(0,n^{1/d})^d$. 
Then for $s>1$,
 arguing as in the proof of Lemma \ref{stab2} yields
\begin{align*} \Pr \left( R_\xi (n^{1/d} \bx ; n^{1/d} \Po_n , n^{1/d} (0,1)^d) > s \right)
\leq \Pr \left( d_* (n^{1/d} \bx ; n^{1/d} \Po_n) > s-1 \right) \\
\leq \exp \left( - \left| A_n^d (\bx, s-1) \right| \right).\end{align*}
So by (\ref{0612b}), for $n \geq n_0$ and $2 < s \leq d^{1/2}n^{1/d} +1$, we obtain,
\[ \sup_{\bx \in \Gamma_n}
\Pr \left( R_\xi (n^{1/d} \bx ; n^{1/d} \Po_n , n^{1/d} (0,1)^d) > s \right)
\leq \exp \left( -C_d (s-1)^{d-1} 
 \right).\]
Also, this probability is zero for $s  > d^{1/2}n^{1/d} +1$.
 Thus
 for any $s > d^{1/2}n_0^{1/d} +1$,
\begin{align*}
\sup_{n \geq 1; ~\bx \in \Gamma_n} 
\Pr \left( R_\xi( n^{1/d} \bx ; n^{1/d} \Po_n , n^{1/d} (0,1)^d) > s \right)
& \leq   \exp \left (-C_d(s-1)^{d-1} \right),\end{align*}
and (\ref{2b}) follows.
$\square$

\section{Boundary effects in the MDST} \label{bdry}

In this section, we consider  the contribution to the total
power-weighted length of the MDST under $\pos$ due to boundary effects
near the `bottom face' of the $d$-cube. Here the possibility of long edges 
leads to rather special behaviour. We shall see that the on-line nearest-neighbour
graph, as described in Section \ref{secong}, will be a useful tool here.

Fix $\eps>0$ small. 
Let $(t_n)_{n > 0}$ be such that $t_n \in (0,1)$ and
$t_n = \Theta (n^{-(1/2)-\eps})$
as $n \to \infty$ (we make a specific choice for $t_n$ in Section \ref{totallength}).
Let $B_n$ denote the boundary region 
$(0,1)^{d-1} \times (0,t_n]$, i.e.~we look
in a thin slice at the base (in the sense of $\pos$) of the unit $d$-cube.
Recall from (\ref{0714a}) that
 $\LL^{d,\alpha} (\X ; R)$ denotes the contribution to the
total weight of the MDST on $\X$ from those points of $\X \cap R$, and
$\tilde \LL^{d,\alpha} (\X; R)
:= 
 \LL^{d,\alpha} (\X; R)
- \Exp [ \LL^{d,\alpha} (\X; R)]$. 
Also recall that
$\Po_n$ denotes a homogeneous Poisson point process of intensity $n$
on $(0,1)^d$. 
Our main result of this section is the following. 

\begin{theorem} 
\label{thmbdry}
Suppose $d \in \{2,3,\ldots\}$.  Let $\eps>0$ and $t_n =
\Theta (n^{-(1/2)-\eps})$ specify $B_n$.
\begin{itemize}
\item[(i)]
Suppose $\alpha \geq d/2$.
With
$Q(d-1,\alpha)$  
as in Lemma \ref{onngthm}, we have that
as $n \to \infty$
\begin{align}
\label{bdry1}
 \tLalph ( \Po_n ; B_n ) \tod Q(d-1,\alpha).
\end{align}
 \item[(ii)] 
 Suppose $\alpha \in (0,d/2)$.
 As $n\to \infty$,
\begin{align} 
\label{0214h}
n^{(\alpha/d)-(1/2)} \tLalph ( \Po_n ; B_n ) \toP 0.
\end{align}
\end{itemize}
\end{theorem}

The idea behind the proof of Theorem \ref{thmbdry} is
to show that the MDST under $\pos$ near to the boundary is close to
an ONG
 defined on a sequence of uniform random vectors in $(0,1)^{d-1}$
coupled to the points of the MDST in $B_n$. To do this,
we produce an explicit sequence of random variables on which we
construct the ONG coupled to
$\Po_n$
on which the MDST is constructed.

Define the point process
\begin{align}
\label{0701f}
 \W_n := \Po_n \cap B_n.
\end{align}
 Let $\beta_n := \card( \W_n )$.
List
$\W_n$ in order of increasing $x_d$-coordinate as
$\bU_i$, $i=1,2,\ldots,\beta_n$. In coordinates, set $\bU_i =
(U_i^1,U_i^2,\ldots,U_i^d)$ for each $i$. Let $\bV_i = (U_i^1,\ldots,U^{d-1}_i)
\in (0,1)^{d-1}$ be the projection of $\bU_i$ down (in the $\be_d$ sense)
onto the base of the unit $d$-cube. 
Set 
\begin{align}
\label{vndef}
\V_n := ( \bV_1, \ldots , \bV_{\beta_n}
).
\end{align} 
Then $\V_n$
is a sequence of uniform random vectors in
$(0,1)^{d-1}$ (the base of the unit $d$-cube), 
on which we may construct the ONG as appropriate. 
Note that the points of $\V_n$ in fact constitute a homogeneous Poisson point process
of intensity $nt_n = \Theta ( n^{(1/2)-\eps})$ on $(0,1)^{d-1}$ (this follows from
the Mapping Theorem, see \cite{kingman}).
With the ONG weight functional
 $\O^{d,\alpha}(\cdot)$ 
defined in Section \ref{secong},
the ONG weight $\O^{d-1,\alpha}(\V_n)$ is coupled in a natural way
to $\LL^{d,\alpha}(\W_n) = \LL^{d,\alpha}(\Po_n;B_n)$. 

Our first step towards Theorem \ref{thmbdry}
is the following result, which shows that, near the boundary, the MDST
is close to the coupled ONG.

\begin{lemma} \label{1108c}
Suppose $d \in \{2,3,\ldots\}$.
Let $\eps>0$ and $t_n =
\Theta (n^{-(1/2)-\eps})$ specify $B_n$.
Let $\W_n, \V_n$ be as defined at
(\ref{0701f}), (\ref{vndef})
respectively.
 For $\alpha \geq 1$, as $n \to \infty$,
\begin{align} 
\label{0209cc} \LL ^{d,\alpha} (\W_n) - \O^{d-1,\alpha}
(\V_n) \to 0, \textrm{ in } L^1,
\end{align} and, 
for $\alpha \in (0,1)$, as $n\to \infty$,
\begin{align} 
\label{0214aa}
\Exp  \left| \LL^{d,\alpha} (\W_n) - \O^{d-1,\alpha} (\V_n ) \right|
 = O \left( n^{(1/2)-\eps-\alpha ((1/2)+\eps)} \right) .\end{align}
\end{lemma} 
\proof
We construct the MDST on the set of points $\W_n$,
 and construct the ONG on the projections down onto $(0,1)^{d-1}$, $\V_n$.
With a slight abuse of notation, 
consider the points $\bU_i=(\bV_i,U^d_i)$, $i=1,\ldots,\beta_n$.

Note that by construction of the MDST 
on $\pos$ and the ONG, and our choice of ordering of points,
we have that
$\bU_j \pos \bU_i$ if and only if 
$j \leq i$. 
Thus
either an edge exists from $\bU_i$ in the MDST
and also from $\bV_i$ in the ONG, or from neither. For the difference between the total weights
of the two models, it suffices to consider the case in which both edges exist.
Then $\bV_i$ is joined to a point
 $\bV_{D(i)}$, $D(i)<i$ in the ONG, and $\bU_i$ to
a point $\bU_{J(i)}$ in the MDST; we do not necessarily have
$J(i)=D(i)$. Since $J(i)<i$ by construction of the MDST on $\pos$ and the ordering
of our points,
we have that $\bV_{J(i)}$ was an admissible candidate to be the point
that $\bV_i$ joins to in the ONG. Therefore, we have that $\| \bV_i-\bV_{D(i)} \|_{d-1} \leq
\| \bV_i-\bV_{J(i)} \|_{d-1}$. It then
follows that
 \begin{align}
 \label{1030y}
 \|(\bV_i,U_i^d)-(\bV_{J(i)},U^d_{J(i)})\|_d^\alpha \geq
\|\bV_i-\bV_{J(i)} \|_{d-1}^\alpha \geq \|\bV_i-\bV_{D(i)} \|_{d-1}^\alpha,
\end{align}
and so we have that, for all $\alpha>0$,
\begin{align}
\label{1002bb}
\O^{d-1,\alpha} ( \V_n ) \leq \LL^{d,\alpha}
( \W_n ) . \end{align}
Also, by construction of the ONG and our ordering on points, we see
 $(\bV_{D(i)},U^d_{D(i)}) \pos (\bV_i,U^d_i)$. So by the construction of the MDST, we have that
\begin{align} 
\label{0209aa}
\| (\bV_i,U^d_i)-(\bV_{J(i)},U^d_{J(i)})
\|_d \leq \| (\bV_i,U^d_i)-(\bV_{D(i)},U^d_{D(i)}) \|_d.
\end{align}
If $(\bx,y) \in (0,1)^{d-1} \times (0,1)$ then $\| (\bx,y) \|_d \leq \|\bx\|_{d-1} + y$, and by the
Mean Value Theorem for the function $t \mapsto t^\alpha$, for $\alpha
\geq 1$,
\begin{align*}
 \| (\bx,y) \|_d^\alpha - \|\bx\|_{d-1}^\alpha \leq (\|\bx\|_{d-1}+y)^\alpha - 
\|\bx\|_{d-1}^\alpha \leq \alpha ((d-1)^{1/2}+1)^{\alpha-1} y ~~~ (\alpha \geq 1). 
\end{align*}
So we have that, for $d \geq 2$ and
$\alpha \geq 1$, there is a finite positive $C$ such that, a.s.,
\begin{align} 
\label{0209bb}
 \| (\bV_i, U^d_i) - (\bV_{D(i)}, U^d_{D(i)}) 
\|_d^\alpha - \| \bV_i - \bV_{D(i)}\|_{d-1}
^\alpha \leq C
(U^d_i - U^d_{D(i)}) .
\end{align}
Then (\ref{0209aa}) and (\ref{0209bb}) yield, for $\alpha \geq 1$, a.s.,
\begin{align}
\label{1002aa}
 \| ( \bV_i, U^d_i ) - (\bV_{J(i)}, U^d_{J(i)} ) \|_d^\alpha
- \|\bV_i - \bV_{D(i)} \|_{d-1}^\alpha
\leq C ( U^d_i - U^d_{D(i)} ) \leq C t_n , \end{align}
which implies that there exist $C, C' \in (0,\infty)$
such that for all $n \geq 1$
\begin{align}
\label{1030t}
\LL^{d,\alpha} (\W_n) -\O^{d-1,\alpha}
(\V_n) \leq C \beta_n  t_n \leq C' \beta_n n^{-(1/2)-\eps} .\end{align}
 Combining (\ref{1002bb}) and (\ref{1030t}) we have that, for $\alpha \geq 1$,
 some $C \in (0,\infty)$ and all $n \geq 1$, a.s.,
\[ \left|   \LL^{d,\alpha}
( \W_n ) - \O^{d-1,\alpha} ( \V_n ) \right| \leq C \beta_n  n^{-(1/2)-\eps}.\]
Taking expectations, using the facts that
$\beta_n$ is Poisson with mean
$nt_n = \Theta (n^{(1/2)-\eps})$ and $\eps >0$,
we obtain (\ref{0209cc}). 

Now we consider the case $\alpha \in (0,1)$. By the concavity of the function $t \mapsto t^\alpha$
for $\alpha \in (0,1)$,
we have for $(\bx,y) \in (0,1)^{d-1} \times (0,1)$ that
 \begin{align*}
 \| (\bx,y) \|_d^\alpha -\|\bx\|_{d-1}^\alpha \leq
(\|\bx\|_{d-1}+y)^\alpha -\|\bx\|_{d-1}^\alpha \leq y^\alpha ~~~ (0 <\alpha <1) .\end{align*}
Then, by a similar argument to the $\alpha \geq 1$ case, we obtain
\[ \left| \LL^{d,\alpha} (\W_n ) - \O^{d-1,\alpha} (\V_n ) \right| \leq C \beta_n n^{-\alpha(
(1/2)+\eps)}, \]
so taking expectations yields (\ref{0214aa}).
 $\square$

\begin{lemma} 
\label{1108d} 
Suppose $d \in \{2,3,\ldots\}$ and $\alpha \geq d/2$.
Let $\W_n$ be as defined at (\ref{0701f})
and suppose that
$Q(d-1,\alpha)$ is the mean-zero
random variable
  given in Lemma \ref{onngthm}. 
Then
as $n \to \infty$
\begin{align}
\label{0715a}
 \tilde \LL^{d,\alpha} (\W_n)  \tod  Q(d-1,\alpha).\end{align}
In particular, $Q(1,\alpha)=\tGal$ for $\alpha \geq 1$, where
$\tGal$ has the distribution given by (\ref{0926a})
or (\ref{0926p}). 
\end{lemma} 
\proof
For $\alpha \geq d/2 \geq 1$, it follows from
(\ref{0209cc}) that as $n \to \infty$
\begin{align}
\label{87b}
\tilde \LL^{d,\alpha} (\W_n) - \tilde \O^{d-1,\alpha} (\V_n) \toP 0.
\end{align}
Also, since $\V_n$ is a homogeneous Poisson point process
of intensity $nt_n = \Theta (n^{(1/2)-\eps})$ on $(0,1)^{d-1}$, and $\alpha \geq d/2 > (d-1)/2$,
we have from 
Lemma \ref{onngthm} that as $n \to \infty$
\begin{align}
\label{87a}
\tilde \O^{d-1,\alpha} (\V_n) \tod Q(d-1,\alpha).\end{align}
Thus (\ref{87b}), (\ref{87a}) and Slutsky's theorem complete the proof of (\ref{0715a}).
$\square$\\

\noindent \textbf{Proof of Theorem \ref{thmbdry}.} For $\alpha \geq d/2$,
(\ref{bdry1}) 
follows from (\ref{0715a}). 
Now suppose $\alpha \in(0,d/2)$.
First suppose that $\alpha <1$. Then (\ref{0214aa}) implies that for $\alpha \in (0,1)$
\begin{align*}
 \Exp  \left[ n^{(\alpha/d)-(1/2)}\left| \LL^{d,\alpha} (\Po_n;B_n) - \O^{d-1,\alpha} (\V_n ) \right|
\right] 
 = O ( n^{\alpha((1/d)-(1/2))-\eps(1+\alpha)} ) \to 0, \end{align*}
since $d \geq 2$ and $\eps>0$. 
So for $\alpha \in (0,1)$ we have
\begin{align}
\label{0022} n^{(\alpha/d)-(1/2)}\left( \tilde \LL^{d,\alpha} (\Po_n;B_n) - \tO^{d-1,\alpha} (\V_n ) \right)
\to 0, \textrm{ in } L^1, \end{align}
as $n \to \infty$.
Also, (\ref{0209cc}) implies that
(\ref{0022}) also holds for $\alpha \in [1,d/2)$ when $d \geq 3$. Thus (\ref{0022}) holds
for all $\alpha \in (0,d/2)$.
Recall that $\V_n$ is a homogeneous Poisson point process in $(0,1)^{d-1}$
with intensity $nt_n = \Theta (n^{(1/2)-\eps})$.
If $\alpha \leq (d-1)/2$, then by Lemma \ref{ongl}(i) and (ii) we have that
for some $C \in (0,\infty)$
\begin{align*} \Var \left[ n^{(\alpha/d)-(1/2)} \tO^{d-1,\alpha} (\V_n ) \right] 
 \leq  C n^{(2\alpha/d)-1} ( n^{(1/2)-\eps} )^{1-(2\alpha/(d-1))} \log n\\
 \leq  C n^{\alpha((2/d)-1/(d-1))-(1/2)} \log n \leq Cn^{-1/d} \log n
\to 0,\end{align*}
as $n \to \infty$.
If $\alpha \in ((d-1)/2,d/2)$, then by Lemma \ref{onngthm}, as $n \to \infty$
\[ n^{(\alpha/d)-(1/2)} \tO^{d-1,\alpha} (\V_n ) \toP 0.\]
So by Slutsky's theorem with (\ref{0022}) we obtain (\ref{0214h}).
$\square$

\section{Proof of Theorem \ref{dickman}}
\label{secdick}

In this section we are
interested in the longest edge in the MDST under $\pos$ on 
$\Po_n \subset (0,1)^d$.
The intuition behind Theorem \ref{dickman} is that
this edge is likely to be near the lower $(d-1)$-dimensional boundary. 
Thus we again make use of the fact that
the MDST near the boundary is well-approximated
by the appropriate on-line nearest-neighbour graph. Then we deduce Theorem \ref{dickman}
from Theorem \ref{ongmax} using the set-up of Section \ref{bdry}.

From Section \ref{bdry} recall that for fixed $\eps>0$,
$B_n$ denotes the boundary region
$(0,1)^{d-1} \times (0, t_n]$ (where $t_n = \Theta ( n^{-(1/2)-\eps})$), and from (\ref{0701f})
that $\W_n = \Po_n \cap B_n$. Also, recall 
from (\ref{vndef})
that $\V_n$ is the sequence
of $(d-1)$-dimensional projections of $\W_n$ in order of increasing $x_d$-coordinate.\\

\noindent{\bf Proof of Theorem \ref{dickman}.}
We have from (\ref{1030y}) 
that every edge in the ONG on $\V_n$ has length
bounded above by the length of some edge in the MDST on $\W_n$. On the other hand,
we have from
(\ref{1002aa})
that  an edge from $\bU_i \in \W_n$ in the MDST
has length at most   $O(t_n)$ greater
 than the edge in the ONG from
the corresponding $\bV_i \in \V_n$. Thus we have that for some $C \in (0,\infty)$ and all $n \geq 1$
\[ 0 \leq \Lm^d (\W_n)  - \Om^{d-1} (\V_n) \leq C n^{-(1/2)-\eps } .\]
Hence almost surely
\[  \Lm^d (\W_n)  - \Om^{d-1} (\V_n) \longrightarrow 0,\]
as $n \to \infty$. By Theorem \ref{ongmax}(i)
and the fact that $\V_n$ is a homogeneous Poisson point process
of intensity $n t_n \to \infty$ (for $\eps$ small), we have
\[ \Om^{d-1} (\V_n) \tod Q_{\rm max} (d-1),\]
as $n \to \infty$. Hence Slutsky's theorem implies that
\begin{align}
\label{hhh1}
  \Lm^d (\W_n)  \tod Q_{\rm max} (d-1),\end{align}
as $n \to \infty$. Set
\[ M_n := \max_{\bx \in \Po_n \setminus \W_n } d_* ( \bx ; \Po_n ),\]
the
length of the longest edge in the MDST
from points of $\Po_n$ in the region
$(0,1)^{d-1} \times (t_n,1)$.
 Then  for any $n \geq 1$,
$\Lm^d (\Po_n) = \max \{ \Lm^d (\W_n) , M_n \}$;
thus
\begin{align}
\label{hhh2}
  \Lm^d (\W_n) \leq \Lm^d (\Po_n) \leq \Lm^d (\W_n) + M_n . \end{align}
Hence by (\ref{hhh2}), (\ref{hhh1}), and Slutsky's theorem,  
to complete the proof of the theorem it suffices
to show that as $n \to \infty$
\begin{align}
\label{hhh3}
M_n \toP 0 .\end{align}
We prove this last statement.
For $\eps >0$ as before
and $(i_1, \ldots, i_d) \in \N^d$, define the cuboid
\[ C(i_1,\ldots,i_d) := ( (i_1-1) \lfloor n^{\eps} \rfloor^{-1} , i_1 \lfloor n^{\eps} \rfloor^{-1} ] \times
\cdots \times ( (i_{d-1}-1) \lfloor n^{\eps} \rfloor^{-1} , i_{d-1}  \lfloor n^{\eps} \rfloor^{-1} ]
\times ( (i_d-1) \lfloor t_n^{-1} \rfloor^{-1}, i_d  \lfloor t_n^{-1} \rfloor^{-1} ] .\]
Let $E_n$ denote the event
\[ \bigcup_{ (i_1,\ldots, i_d) \in \N^d \cap [ (0, \lfloor n^\eps\rfloor]^{d-1}\times (0, \lfloor t_n^{-1} \rfloor] ]}
\left\{ \Po_n \cap  C(i_1,\ldots, i_d) = \emptyset \right\}.\]
The number of points of $\Po_n$ in each cuboid $C(i_1,\ldots,i_d)$ in the 
union
is Poisson distributed
with mean 
\[ n \cdot \lfloor n^\eps \rfloor^{-(d-1)} \cdot
 \lfloor t_n^{-1} \rfloor^{-1}   = \Theta ( n^{(1/2)-d\eps}), \]
and the total number of cuboids in the union is 
$\lfloor n^\eps \rfloor^{d-1} \lfloor t_n^{-1} \rfloor = O( n^{(1/2) +d \eps} )$. 
 Thus Boole's inequality implies that there exist $C, C' \in (0,\infty)$
 for which, for all $n \geq 1$,
\[ \Pr (E_n) \leq C n ^{(1/2)-d\eps} \exp ( - C' n^{(1/2)-d\eps} ) ,\]
and hence $\Pr(E_n) \to 0$ 
as $n \to \infty$, for $\eps$ small enough. However, if $E_n$ does not occur then each
cuboid contains at least one point of $\Po_n$ and $M_n$ is bounded by a constant
times $n^{-\eps}$. Thus (\ref{hhh3}) follows and the proof is complete. $\square$

\section{Proof of Theorem \ref{mainth}} \label{totallength}

In this section we complete the proof of our convergence
in distribution result for $\LL^{d,\alpha}(\Po_n)$, Theorem \ref{mainth}. 
Recall from Section \ref{middle} that $\eps>0$ is fixed (small)
and $\Gamma_n$ denotes the region
 $(0,1)^{d-1} \times (g_n,1)$, where  
$g_n =\Theta(n^{\eps-(1/d)})$ as $n \to \infty$.
As in Section
\ref{bdry}, denote
by $B_n$
the region $(0,1)^{d-1} \times (0, t_n]$, where $t_n = \Theta (n^{-(1/2)-\eps})$. 
 We will make a particular choice for $g_n$ and $t_n$ shortly.
 Denote by $I_n$ the intermediate region $(0,1)^d \setminus ( B_n \cup \Gamma_n) = (0,1)^{d-1} \times 
 (t_n,
g_n]$.

In order to prove Theorem \ref{mainth}, we need to collect previous results on the 
limiting behaviour of the MDST in the regions $\Gamma_n$ and $B_n$, and also
deal with the region $I_n$.
In Sections \ref{ltot} and \ref{bdry} we saw
 that, for large $n$, the weight
 (suitably centred and scaled)
  of
edges starting in $\Gamma_n$ satisfies
 a central limit theorem, and the weight
of edges starting in
 $B_n$ can be approximated by the on-line nearest-neighbour graph.
To complete the proof of 
Theorem \ref{mainth}, we shall
  show that (with a suitable
scaling factor for $\alpha<d/2$) the contribution
to the total weight from points in $I_n$ 
has variance converging to 
 zero, and that the lengths from $B_n$ and $\Gamma_n$
are asymptotically independent by virtue of the fact that
the configuration of points in $I_n$ is (with probability
approaching one) sufficient to ensure that the configuration
of points in $B_n$ has no effect on the edges from points in $\Gamma_n$.

Recall from (\ref{0714a})
 that for a point set $\X \subset \R^d$
and a region $R \subseteq \R^d$,
 $\LL^{d,\alpha}(\X;R)$ denotes
the total weight of edges of the MDST on $\X$ which originate
in the region $R$.
The next result is the main result of this section: it gives asymptotic
control of the 
variance of $\LL^{d,\alpha} (\Po_n;I_n)$, and will allow us to complete the proof of Theorem
\ref{mainth}.
\begin{lemma}
\label{99c}
Suppose $d \in \{2,3,4,\ldots\}$ and $\alpha >0$. Then for
small enough $\eps>0$ there exist
$g_n = \Theta (n^{\eps-(1/d)})$ and $t_n=\Theta(n^{-(1/2)-\eps})$
specifying $I_n$
for which,
as $n \to \infty$,
\begin{align}
\label{99ca}
 \Var [ {\cal L}^{d,\alpha} (\Po_n ; I_n ) ] \to 0,~~~(\alpha > (d-1)/2 ), \\
\label{99cb}
\textrm{ and } ~~~  \Var [ n^{(\alpha/d)-(1/2)} {\cal L}^{d,\alpha} (\Po_n ; I_n ) ] \to 0,~~~(0 < \alpha < d/2).\end{align}
\end{lemma}

Before embarking on the proof of Lemma \ref{99c}, we prove the following preliminary result which,
for our purposes, will control the dependency structure of the MDST.
Let $\X$ be a set of points in $(0,1)^d$. For 
non-empty $\X$ and
$\bx \in \X$,
let $D_*(\bx;\X)$ denote the total
degree  of $\bx$ 
(i.e.~the total number of
directed edges that have $\bx$ as one endpoint)
in the MDST on $\X$; set
$D_*(\bx ; \emptyset) :=0$ for any $\bx$.

\begin{lemma}
\label{degs}
Let $d \geq 2$. For any $\eps \in (0,1)$
there exist $C, C' \in (0,\infty)$ such that
for all $n \geq 1$
\[ \Pr \left( \sup_{\bx \in \Po_n} D_*(\bx ; \Po_n) > n^\eps \right) \leq C \exp( - C' n^\eps).\]
\end{lemma}
\proof
Suppose $d \geq 2$. 
Fix $n \in \N$. Let $(\bU_1,\ldots,\bU_n)$
be the points of a binomial point process
of $n$ independent uniform random vectors on $(0,1)^d$ ,
listed in order of increasing $x_d$-coordinate,
so that $\bU_1 \pos \bU_2 \pos \cdots \pos \bU_n$. 
Set $\X_n := \{ \bU_1,\ldots, \bU_n\}$.
 
We now consider our usual
coupling of the MDST to the ONG. 
In coordinates, write $\bU_i=(U_i^1,\ldots,U_i^d)$. Set $\bV_i = (U_i^1,\ldots,U_i^{d-1})$, the
projection of $\bU_i$ down (in the $\be_d$-sense) onto $(0,1)^{d-1}$. 

With probability one, the $\bU_j$, $\bV_j$ have distinct
$d$-, $(d-1)$-dimensional inter-point distances, so
there are no ties to break in constructing the MDST
or ONG.
Consider a point $\bU_j$ with $j \in \{1, \ldots, n - 1\}$.
Suppose that $\bU_k$, $j < k \leq n$ is joined
to $\bU_j$ in the MDST on $\X_n$. Then
$\| \bU_k - \bU_j \|_d \leq \| \bU_k -\bU_i \|_d$
for $i \in \{ j+1,\ldots, k-1\}$. Also
\[ \| \bV_k - \bV_i \|_{d-1}^2 =  \| \bU_k - \bU_i \|_d^2
-( U_k^d - U_i^d)^2 .\]
Then since $U_i^d$ is increasing in $i$,
 $\| \bV_k - \bV_i \|_{d-1}$ is minimized
over $i \in \{j,\ldots,k-1\}$ by $i=j$.
In other words,  
a necessary condition
for $\bU_k$, $j < k \leq n$, to be joined to $\bU_j$ 
in the MDST on $\X_n$ is that the corresponding
edge from $\bV_k$ to $\bV_j$ 
exists in the ONG on sequence of points $(\bV_j,\bV_{j+1},\ldots,\bV_n)$
in $(0,1)^{d-1}$. 

Hence the in-degree of
$\bU_j$ in the MDST on $\X_n$ is bounded above
by the in-degree of $\bV_j$ in the ONG on $(\bV_j,\bV_{j+1},\ldots,\bV_n)$.
Since $\bV_1,\ldots,\bV_n$ are independent uniform random vectors in $(0,1)^{d-1}$,
we have that this latter quantity has the same distribution as 
 the degree of $\bV_1$ in the ONG on $(\bV_1,\bV_2,\ldots,\bV_{n-j+1})$.
 Hence $D_* (\bU_j ; \X_n)$ is stochastically dominated by
 the degree of $\bV_1$ in the ONG on $(\bV_1,\bV_2,\ldots,\bV_{n})$,
 which we denote
 $D_{\rm ONG} (n)$.

Hence
\[ \sup_{1 \leq j \leq n} \Pr (D_*(\bU_j;\X_n )>s) \leq \Pr ( D_{\rm ONG} (n) >s).\]
Then by Boole's inequality, we have that
\begin{align*} \Pr \left( \sup_{1 \leq j \leq n} D_*(\bU_j ; \X_n) > s \right) 
 \leq  \sum_{j=1}^n \sup_{1 \leq i \leq n} \Pr (  D_*(\bU_i ; \X_n) > s )
 \leq  n \Pr ( D_{\rm ONG} (n) >s).\end{align*}
Let $N(n)=\card(\Po_n)$. We have
\begin{align*}
\Pr \left( \sup_{1 \leq j \leq N(n)} D_*(\bU_j ; \X_{N(n)}) > s \right)
& \leq   \Pr( N(n) \geq 2n) + \sup_{m < 2n} \Pr \left( \sup_{1 \leq j \leq m} D_*(\bU_j ; \X_m) > s \right)\\
& \leq   \Pr( N(n) \geq 2n) +  2n \Pr ( D_{\rm ONG} (2n) >s).
\end{align*}
By  following the
argument in Section 3.1 of \cite{boll}
we have that for any $\eps>0$,
$\Pr ( D_{\rm ONG} (2n) > n^\eps ) = O( \exp ( -Cn^\eps))$. 
Also, $\Pr (N(n) \geq 2n) = O( \exp ( -Cn ) )$ by standard Poisson tail bounds
(e.g.~Lemma 1.2 in \cite{penbook}). 
This completes the proof.
$\square$ 
\\

To prove Lemma \ref{99c} we first derive an upper bound
(\ref{vv1} below) for $\Var [ \LL^{d,\alpha} (\Po_n ; I_n)]$
in terms of the mean-square changes in $\LL^{d,\alpha} (\Po_n ; I_n)$
on re-sampling Poisson points over a certain partition of
$B_n \cup I_n$ into boxes, in a similar way to
a technique in \cite{penyuk1}.
 Unlike in \cite{penyuk1},
where the boxes are the same shape and size, we need to
use boxes of different shapes to take account of the structure of the MDST near the boundary.
 
For each $n \geq 1$, we will
 divide $(0,1)^d$ into layers of rectangular $d$-cells. 
 To begin we will divide $(0,1)^{d-1} \times (0,\infty)$ into layers
  starting at the base (in the $\be_d$ sense).
  The $k$-th layer $(k \in \N)$ will have
 height $h_n(k)$ given by
\[ h_n(k) :=  
n^{-1+\eps} 2^{(k-1)(d-1)}. \]
We will let
$H_n(k)$ denote the starting height (in the $\be_d$ sense)
of layer $k$; define
$H_n(1):=0$ and for $k \geq 2$
\[ H_n(k) := \sum_{i=1}^{k-1} h_n(i)
=  \sum_{i=0}^{k-2} n^{-1+\eps} 2^{(d-1)i} 
= c_d n^{-1+\eps} \left( 2^{(d-1)(k-1)}  - 1 \right)
= c_d h_n (k) - c_d n^{-1+\eps},\]
where $c_d = (2^{d-1}-1)^{-1}$ depends only on $d$.
We then define the box
\[ L_n(k) := (0,1)^{d-1} \times (H_n(k),H_n(k+1)] ~~~(k \in \N);\]
we will refer to $L_n(k)$ as the $k$-th {\em layer}. 

For $n \geq 1$ define $M_n \in \N$ such that 
\[ M_n := \min \{ m \in \N : H_n(m+1) \geq n^{-(1/2)-\eps} \}.\]
Then $M_n$ satisfies
\begin{align}
\label{mmm}
 M_n = \Theta (\log n), ~~~ 2^{M_n} = \Theta (n^{(1-4\eps)/(2(d-1)) }).\end{align}
We then define  for $n \geq 1$ the region
\[ B_n := \bigcup_{k=1}^{M_n} L_n(k) = (0,1)^{d-1} \times (0, H_n(M_n +1)] .\]
Then with our previous notation as $B_n = (0,1)^{d-1} \times (0, t_n)$,
we have  $t_n = H_n(M_n +1) = \Theta ( n^{-(1/2)-\eps} )$.

Also for $n \geq 1$
define $K_n \in \N$ such that
\[ K_n := \min \{ k \in \N : H_n(k+1) \geq n^{\eps-(1/d)}\}.\]
Thus 
\begin{align}
\label{kkk}
 K_n  := \Theta ( \log n) , ~~~ 2^{K_n} = \Theta( n^{1/d}).\end{align}
Define for $n \geq 1$ the region
\begin{align}
\label{indef}
 I_n := \bigcup_{k=M_n+1}^{K_n} L_n(k) = (0,1)^{d-1} \times (H_n(M_n+1),H_n(K_n+1)],\end{align}
so that, with our previous
notation for $I_n$,
 $g_n = H_n(K_n+1) = \Theta (n^{\eps -(1/d)})$. These specific choices
for $t_n$ and $g_n$ then fit with our previous usage.

We now subdivide each layer into cells.
For $k=1,2,\ldots,K_n$, divide layer $k$ into rectangular cells
of height $h_n(k)$ by forming a grid by
dividing each of the $d-1$ sides of the layer into $2^{k-1}$
equal intervals. Thus layer $k$ then consists of $2^{(k-1)(d-1)}$
cells of height $h_n(k)$
and $(d-1)$-widths $2^{1-k}$. Each such cell has volume
$2^{(1-k)(d-1)}h_n(k)=n^{-1+\eps}$. The total number of cells in all the layers up to layer $K_n$ we denote by
$\ell(n)$, which is given by
\begin{align}
\label{555}
 \ell(n) := \sum_{k=1}^{K_n}  2^{(k-1)(d-1)}  = 
 \Theta(2^{(d-1)K_n}) = \Theta(n^{1-(1/d)}),\end{align}
 by (\ref{kkk}).
For layers
$1$ up to $K_n$, label the individual cells
lexicographically as $S^n_i$, $i=1,2,\ldots,\ell(n)$. 

Note that for $\eps$ small enough, cells in layer $k \leq M_n$ are always
wider than they are tall, while for $M_n \leq k \leq K_n$ cells
in layer $k$ have height at most a constant times $n^\eps$ times
their width.

Let
$\tilde \Po_n$ denote an independent copy of the homogeneous Poisson point 
process $\Po_n$, and for $i=1,2,\ldots,\ell(n)$ set
\[ \Po_n^i := ( \Po_n \setminus S^n_i ) \cup ( \tilde \Po_n \cap S^n_i ),\]
so that $\Po_n^i$ is $\Po_n$ but with the Poisson points in
 $S^n_i$ independently
re-sampled. For ease of notation during this proof, for $n >0$ set 
$Y_n = \tilde \LL^{d,\alpha} (\Po_n; I_n)$. Define
\[ \Delta^n_i := \tilde \LL^{d,\alpha} ( \Po^i_n ; I_n ) - \tilde \LL^{d,\alpha}
(\Po_n ; I_n) = \LL^{d,\alpha} ( \Po^i_n ; I_n ) - \LL^{d,\alpha}
(\Po_n ; I_n),\]
the change in $Y_n$ on re-sampling the Poisson points in $S^n_i$. 
 By Steele's \cite{steeleES} version of the Efron--Stein
inequality, or by a martingale difference argument,
we have that for $n >0$
\begin{align}
\label{vv1}
 \Var [ \LL^{d,\alpha} (\Po_n ; I_n ) ] =
 \Exp[ Y_n^2] \leq
  \sum_{i=1}^{\ell(n)} \Exp[(\Delta^n_i)^2].\end{align}
For $i=1,2,\ldots, \ell(n)$, let $G(i)$ be the integer $k \in \{1, \ldots, K_n\}$
 such that $S_i \subseteq L(k)$, so that $G(i)$ is the layer
to which $S_i$ belongs. Formally,
\[ G(i) := \left\lceil (d-1)^{-1} \log_2 \left[ \left( 2^{d-1} -1\right) i + 1 \right] \right\rceil .\]
The next result gives bounds on
$\Exp[ (\Delta_i^n)^2 ]$.

\begin{lemma}
\label{lem1008}
Let $d \in \{2,3,\ldots\}$ and $\alpha >0$. There exists $C \in (0,\infty)$
such that for all $n>0$ and all $i \in \{1,2,\ldots,\ell(n)\}$
\begin{align}
\label{1008a}
\Exp [ (\Delta_i^n)^2 ] \leq \begin{cases}
C n^{(6+4\alpha) \eps} n^{-\alpha/(d-1)} & \textrm{~if~} G(i) \leq M_n \\
C n^{(6+2\alpha)\eps} 2^{-2\alpha G(i)} & \textrm{~if~} M_n < G(i) \leq K_n
\end{cases}
\end{align}
\end{lemma}
Note that $2^{-2\alpha G(i)} = \Theta (i^{-2\alpha/(d-1)})$ as
$i \to \infty$, and for $G(i) \leq M_n$ or $G(i) \leq K_n$,
$ i = O( n^{(1-4\eps)/2} )$ or $i = O (n^{1-(1/d)})$ respectively.\\

\noindent 
{\bf Proof of Lemma \ref{lem1008}.}
Let $E'_n$ denote the event that every cell $S^n_j \subset (B_n \cup I_n)$
contains at least one and not more than $n^{2\eps}$ points of $\Po_n$, and
also $\tilde \Po_n$.
That is,
\begin{align*}
 E'_n := \bigcap_{1 \leq j \leq \ell(n)} 
 \left\{ 1 \leq \card ( \Po_n \cap S^n_j) \leq n^{2\eps},
 ~ 1 \leq \card ( \tilde \Po_n \cap S^n_j ) \leq n^{2 \eps} \right\}.\end{align*}
 We have  from Boole's inequality and the fact that
  $\card (\Po_n \cap S^n_j)$ has the same distribution
  as  $\card (\tilde \Po_n \cap S^n_j)$ 
  \begin{align}
  \label{1008b}
 \Pr ((E'_n)^c) & \leq   
2 \sum_{1 \leq j \leq \ell(n) }
\Pr \left( \left\{ 1 \leq \card ( \Po_n \cap S^n_j) \leq n^{2\eps} \right\}^c \right) \nonumber\\
& = 2 \ell (n) \left[ \Pr ( \card ( \Po_n \cap S^n_j) > n^{2\eps} )
+ \Pr ( \card ( \Po_n \cap S^n_j) = 0 ) \right] .\end{align}
Now $\card (\Po_n \cap S^n_j)$, $j = 1, \ldots, \ell(n)$
are Poisson distributed with mean $n^{\eps}$ (since
$|S_j| = n^{-1+\eps}$).
 By standard Chernoff bounds on Poisson tails (see e.g.~Lemma 1.2
of \cite{penbook}), we have that $\Pr(\card (\Po_n \cap S^n_j) > n^{2\eps}) = O( \exp (-Cn^{2\eps} \log n))$,
whereas $\Pr(\card (\Po_n \cap S^n_j) =0)=   \exp (-  n^{\eps})$. Thus from (\ref{1008b}), using (\ref{555}),
 there exists $C  \in (0,\infty)$ such that
\begin{align}
\label{1008c}
 \Pr ((E'_n)^c) = O \left( n^{1-(1/d)} \exp(-   n^\eps) \right)  
 =  O (\exp(- C n^\eps)),\end{align}
as $n \to \infty$. 

Also, for $\eps>0$ and $n>0$
 let $E_n''$ denote the event that the maximum vertex degree in the MDST on $\Po_n$ 
and on $\Po_n^i$ for each $i$ is bounded by $n^\eps$; i.e.
\[ E_n'' := \left\{ \sup_{ \X \in \{ \Po_n, \Po_n^1,
\ldots, \Po_n^{\ell (n)} \} }
\sup_{\bx \in \X} D_* (\bx ; \X ) \leq n^{\eps} \right\}.\]
 Then by Lemma \ref{degs} we have that 
 for some $C \in (0,\infty)$,
 \begin{align}
 \label{1008d}
 \Pr ((E_n'')^c) = O(\exp(-Cn^\eps)).
 \end{align} 
Let 
\begin{align}
\label{bxbx}
E_n := E_n' \cap E_n''.\end{align}
 Then $\Pr (E_n^c) \leq \Pr ((E_n')^c )+ \Pr ((E_n'')^c )$ so that
 by (\ref{1008c}) and (\ref{1008d}) we have
 that there exists $C \in (0,\infty)$ such that
 as $n \to \infty$
 \begin{align}
 \label{1008e}
 \Pr (E_n^c) = O(\exp(-C n^\eps)).
 \end{align}

We bound $\Exp [ (\Delta_i^n )^2 ]$ by partitioning over the occurrence  of $E_n$
and using the fact that
\begin{align}
\label{vv22}
 \Exp[ (\Delta^n_i)^2 ] \leq \Exp[ (\Delta^n_i)^2 \mid E_n] + \Exp[ (\Delta^n_i)^2 \1_{E_n^c} ] 
  .\end{align}
First note that
  by the Cauchy--Schwarz inequality
  and the trivial bound
  $| \Delta_i^n| \leq C (\card (\Po_n) + \card (\tilde \Po_n))$,
  we have that
  \[ \Exp[ (\Delta^n_i)^2 \1_{E_n^c} ]
  \leq ( \Exp [ (\Delta_i^n)^4 ] )^{1/2} (\Pr (E_n^c ) )^{1/2}
\leq C ( \Exp [ (N(n)+  N'(n))^4] )^{1/2} (\Pr (E_n^c ) )^{1/2} ,\]
where $N(n), N'(n)$ are independent Poisson random
variables with mean $n$.  
Hence from (\ref{1008e}) we have that
for some $C \in (0,\infty)$
\begin{align}
\label{1008f}
\Exp[ (\Delta^n_i)^2 \1_{E_n^c} ] = O(   \exp (-C n^{\eps}) ).
\end{align}

Next we treat the case where $E_n$
occurs. First suppose $G(i) \leq M_n$, so that $S^n_i \subseteq B_n$. Contributions to $\Delta^n_i$ are 
from directed edges from Poisson points in $I_n$ to Poisson points
in $S_i^n$: specifically, such edges
 that are added or deleted on the re-sampling of the Poisson points
in $S^n_i$. The number of such edges is bounded by the sums
of the vertex degrees in the MDST of points of $\Po_n \cap S_i^n$
and $\tilde \Po_n \cap S_i^n$. Given 
 $E_n$, the number of points of $\Po_n \cap S^n_i$ is bounded by $n^{2\eps}$,
similarly with $\tilde \Po_n$, and each point has degree bounded by $n^{\eps}$. It follows
that the number of edges that can contribute to  $\Delta^n_i$
is bounded by $2n^{3\eps}$ under $E_n$.
Further,
given $E_n$, the length of an edge   contributing
 to $\Delta^n_i$ is bounded by a constant times the width
of cells in $L(M_n+1)$ the first layer in $I_n$, 
which for $d \geq 2$
is $O(2^{-M_n})=O(n^{2\eps-(1/(2(d-1)))})$ by (\ref{mmm}). Each edge therefore
gives a contribution to $\Delta^n_i$
at most $O(n^{2\alpha \eps-(\alpha/(2(d-1)))})$ in absolute value.
It follows that there exists $C \in (0,\infty)$ such that
for all $n>0$ and all $i$ with $G(i) \leq M_n$
\begin{align}
\label{1008g}
 \Exp[ (\Delta^n_i)^2 \mid E_n] \leq C n^{(6+4\alpha)\eps} n^{-\alpha/(d-1)}.\end{align}
Thus from (\ref{vv22}) with the bounds
(\ref{1008f}) and (\ref{1008g}) we obtain
the $G(i) \leq M_n$ case of (\ref{1008a}).
 
 Finally suppose $M_n < G(i) \leq K_n$, so that $S^n_i \subseteq I_n$.
 Given $E_n$, the number of points of $\Po_n \cap S^n_i$ is bounded by $n^{2\eps}$;
 similarly for $\tilde \Po_n$. Further,
given $E_n$, edge lengths contributing to $\Delta^n_i$ are bounded by a constant times
$n^\eps$ times the width
of cell $S_i^n$ in layer $G(i)$, which is $O(2^{-G(i)})$,
and each point has degree bounded by $n^{\eps}$. 
Thus for $M_n < G(i) \leq K_n$,
\begin{align}
\label{vv3}
 \Exp[(\Delta_i^n)^2 \mid E_n ] = O (n^{(6+2\alpha)\eps} \cdot 2^{-2\alpha G(i)}).\end{align}
 Then (\ref{vv22}) with
(\ref{1008f}) and (\ref{vv3}) yields the
$M_n < G(i) \leq K_n$ case of
(\ref{1008a}). $\square$\\

We can now complete the proof of Lemma \ref{99c}.\\

\noindent
{\bf Proof of Lemma \ref{99c}.}
Fix $d \geq 2$ and $\alpha>0$.
Take $I_n$ as defined by (\ref{indef}) so that
$g_n = H_n (K_n+1)$ and
$t_n = H_n (M_n+1)$ are as in
the statement of Lemma \ref{99c}. Again using
the shorthand $Y_n = \tilde \LL^{d,\alpha} (\Po_n ; I_n)$,
we obtain from (\ref{vv1}) with 
(\ref{1008a})
 that for all $n >0$
 \begin{align*}
  \Exp[Y_n^2] & =   \sum_{1 \leq i \leq \ell(n)} \Exp[(D^n_i)^2 ]
  \leq \sum_{k=1}^{M_n}  \sum_{i:S^n_i \subseteq L(k)} \Exp[(\Delta^n_i)^2 ] + 
  \sum_{k=M_n+1}^{K_n}
 \sum_{i:S^n_i \subseteq L(k)} \Exp[(\Delta^n_i)^2 ] \\
 & \leq   C \sum_{k=1}^{M_n} 2^{k(d-1)} n^{(6+4\alpha)\eps} n^{-\alpha/(d-1)}
 + C \sum_{k=M_n+1}^{K_n} 2^{k(d-1)} n^{(6+2\alpha)\eps} 2^{-2\alpha k} \\
 & \leq   C 2^{M_n (d-1)} n^{(6+4\alpha)\eps}n^{-\alpha/(d-1)} + C 2^{K_n (d-1-2\alpha)}n^{(7+2\alpha)\eps}
+ C 2^{M_n (d-1-2\alpha)}n^{(7+2\alpha)\eps}
 ,\end{align*}
 where the additional $n^\eps$ factor in the last two terms
takes care of the extra logarithmic factor
when $\alpha = (d-1)/2$. 
Using (\ref{mmm}) and (\ref{kkk}) 
we thus
have that for any $\eps>0$
 there exists $C \in (0,\infty)$
such that for all $n >0$ 
\begin{align}
\label{ddff}
\Exp[Y_n^2] \leq  
C n^{(1/2)-(\alpha/(d-1)) + (4+4\alpha) \eps} (1+ n^{(1+2 \alpha) \eps}  ) 
 + C n^{1-((1+2\alpha)/d)+(7+2\alpha)\eps}
 .\end{align} 
 For $d \geq 2$,
this tends to zero as $n \to \infty$
 for $\alpha > (d-1)/2$ and $\eps$ sufficiently small, which gives (\ref{99ca}).
On the other hand, for $\alpha < d/2$, we have from (\ref{ddff}),
 noting that $(2\alpha/d)-(\alpha/(d-1)) = (\alpha/d)(d-2)/(d-1)$, that
 \begin{align*}
  \Exp[n^{(2\alpha/d)-1} Y_n^2] 
  \leq  
  C n^{(\alpha/d)(d-2)/(d-1) -(1/2) + (4+4\alpha)\eps} (1 + n^{(1+2\alpha) \eps}  ) 
   + C n^{-(1/d)+(7+2\alpha) \eps}
 ,\end{align*}
 which also tends to zero as $n \to \infty$
 for $\eps$ small enough and $d \geq 2$. 
This gives (\ref{99cb}).
$\square$\\

\noindent
{\bf Proof of Theorem \ref{mainth}.} 
Again we use the construction of Lemma \ref{99c}.
For the duration of this proof, to ease notation, set $X_n = \tilde \LL^{d,\alpha}
(\Po_n;\Gamma_n)$, $Y_n = \tilde \LL^{d,\alpha}
(\Po_n;I_n)$ and $Z_n = \tilde \LL^{d,\alpha}
(\Po_n;B_n)$. 
Thus $\tilde \LL^{d,\alpha}
(\Po_n) = X_n +Y_n +Z_n$.

First suppose $\alpha \in (0,d/2)$. Then from (\ref{0214h}) 
and (\ref{99cb})
we have that $n^{(\alpha/d)-(1/2)} (Y_n + Z_n) \toP 0$ as $n \to \infty$. With 
Lemma
\ref{CLT} and Slutsky's theorem, we obtain the $\alpha \in (0,d/2)$ case of
Theorem \ref{mainth}.

Now suppose $\alpha >d/2$. Then Lemma \ref{CLT} and (\ref{99ca}) imply that
$X_n+Y_n \toP 0$ as $n \to \infty$. So (\ref{bdry1}) with Slutsky's theorem 
gives the $\alpha>d/2$ case of Theorem \ref{mainth}.

Finally suppose $\alpha=d/2$. Again
(\ref{99ca}) implies that $Y_n \toP 0$. Here we have  from (\ref{bdry1}) that $Z_n \tod Q(d-1,d/2)$ and
from Lemma \ref{CLT} that $X_n \tod W_1$ where $W_1$ is a normal random variable.
We need to show that the limits $W_1$ and $Q(d-1,d/2)$ are independent.
Set $k_n := \lceil n^{(1/d)-(\eps/2)} \rceil$. For $\bz \in \Z^{d-1} \cap
[0,k_n]^{d-1}$ let $C(\bz) \subset I_n$ denote
the cube
\[ C(\bz) := (k_n^{-1}\bz,0) + (0,k_n^{-1}]^{d-1} \times (g_n - k_n^{-1},g_n].\]
Thus there are $k_n^{d-1} =\Theta( n^{1-(1/d)-\eps(d-1)/2})$ such cubes,
each of volume $k_n^{-d} =\Theta ( n^{-1+(d\eps/2)})$. Let $A_n$ denote the event
\[ A_n :=  \bigcap \left\{ \card( \Po_n \cap C(\bz) ) > 0 : \bz \in \Z^{d-1} \cap
[0,k_n]^{d-1} \right\} .\] The number of points of $\Po_n$ in each cube
$C(\bz)$ is a Poisson random variable with mean $\Theta(n^{d\eps/2})$, and so
\[ \Pr (A_n^c) \leq \sum_\bz \Pr ( \card( \Po_n \cap C(\bz) ) = 0 ) = O (
n^{1-(1/d)-\eps(d-1)/2} \cdot \exp(-C n^{d\eps/2})) \to 0,\]
as $n \to \infty$. Given a configuration of $\Po_n$ satisfying $A_n$, for $n$ sufficiently large, $X_n$ and $Z_n$ are (conditionally)
independent, since no point of $\Po_n \cap \Gamma_n$ can be joined to a point
of $\Po_n \cap B_n$ in the MDST.
Then the proof is completed by following
the argument for Equation
(7.25) in \cite{total}.
 $\square$

\section{Proof of Theorem \ref{llnthm}}
\label{prflln}

In order to complete the proof of Theorem \ref{llnthm}, we need
to add to the law of large numbers away from the boundary (in region
$\Gamma_n$), Lemma 
\ref{llnlem}, by dealing with the edges near to the boundary. We proceed in a 
similar fashion to Sections \ref{bdry} and \ref{totallength}, dealing 
with the contributions from the region $B_n$ in Lemma \ref{kk1} below
(using the coupling to the ONG as in Section \ref{bdry}),
and with the contributions from the region $I_n$ in Lemma \ref{expClem} below 
(using the construction of Section \ref{totallength}).

\begin{lemma}
\label{expClem}
Suppose $d \in \{2,3,\ldots\}$
and $\alpha > 0$. Then for
small enough $\eps>0$ there exist
$g_n = \Theta (n^{\eps-(1/d)})$ and $t_n=\Theta(n^{-(1/2)-\eps})$
specifying $I_n$
for which,
as $n \to \infty$, 
\begin{align}
\label{101}
n^{(\alpha/d)-1} \LL^{d,\alpha} (\Po_n;I_n) \to 0, \textrm{ in } L^1, ~~~(\alpha \in (0,d)),\\
\label{1111d}
\textrm{ and }~~~
 \LL^{d,\alpha} ( \Po_n ; I_n)  \to 0, \textrm{ in } L^1, ~~~(\alpha > d-1). 
\end{align}
\end{lemma}
\proof 
Recall the construction of the partition
of $I_n$ described in Section \ref{totallength},
and the definition of the event $E_n$ from
(\ref{bxbx}).
Then
\begin{align}
\label{xx0}
 \Exp [ \LL^{d,\alpha} ( \Po_n ; I_n ) ] = \Exp [ \LL^{d,\alpha} ( \Po_n ; I_n ) \1_{E_n} ]
+ \Exp [ \LL^{d,\alpha} ( \Po_n ; I_n ) \1_{E^c_n} ],\end{align}
where by Cauchy--Schwarz
\[ \Exp [ \LL^{d,\alpha} ( \Po_n ; I_n ) \1_{E^c_n} ] \leq 
(\Exp [ ( \LL^{d,\alpha} ( \Po_n ; I_n ) )^2 ])^{1/2} ( \Pr (E^c_n) )^{1/2} 
\leq C (\Exp [ N(n)^2 ])^{1/2} ( \Pr (E^c_n) )^{1/2} , \]
where $N(n) = \card ( \Po_n )$
is Poisson distributed with mean $n$. Thus by
(\ref{1008e}) there exists $C \in (0,\infty)$ such that
\begin{align}
\label{xx1}
  \Exp [ \LL^{d,\alpha} ( \Po_n ; I_n ) \1_{E^c_n} ] = O ( \exp ( -C n^{\eps} ) ).\end{align}
Also, using the construction of Section \ref{totallength},
\begin{align*} \Exp [ \LL^{d,\alpha} ( \Po_n ; I_n ) \1_{E_n} ]
\leq
\sum_{k=M_n+1}^{K_n} \sum_{i:S^n_i \subseteq L(k)} \Exp [ \LL^{d,\alpha} (\Po_n ; S_i^n ) \mid E_n ]
.\end{align*}
Given $E_n$, as in the proof of Lemma \ref{lem1008}, the number of points in each $S_i^n$
is bounded by $n^{2\eps}$, the degree of each point
is bounded by $n^\eps$, and each edge has length bounded by
a constant times $n^\eps 2^{-G(i)}$. Thus 
\begin{align}
\label{xx2}
 \Exp [ \LL^{d,\alpha} (\Po_n ; I_n ) \mid E_n ] 
\leq C
\sum_{k=M_n+1}^{K_n} 2^{k(d-1)} \cdot n^{(3+\alpha) \eps} \cdot 2^{-\alpha k}. \end{align}
Thus from (\ref{xx0}) with (\ref{xx1}) and (\ref{xx2})
we obtain
\[ \Exp [ \LL^{d,\alpha} (\Po_n ; I_n ) ] = O ( 2^{(d-1-\alpha)K_n}  n^{(3+\alpha)\eps})
+ O(2^{(d-1-\alpha)M_n}  n^{(4+\alpha)\eps}),\]
where the additional $n^\eps$ factor in the second term
takes care of the extra logarithmic factor
when $\alpha = d-1$. Using (\ref{mmm}) and (\ref{kkk}) we have
for $d \geq 2$
\begin{align}
\label{xx3}
 \Exp [ \LL^{d,\alpha} (\Po_n ; I_n ) ] = 
O(n^{1-(\alpha/d)-(1/d)+(3+\alpha)\eps})
+O(n^{(1/2)-(\alpha/(2(d-1)))+(2+6\alpha)\eps}).\end{align}
For $\alpha > d-1$ this
tends to zero as $n \to \infty$
 for $\eps$ small enough, and so we obtain
(\ref{1111d}). On the other hand,
for $\alpha \in (0,d)$, we have from (\ref{xx3}) that
\[ \Exp [ n^{(\alpha/d)-1} \LL^{d,\alpha} (\Po_n ; I_n ) ] 
= O( n^{(3+\alpha)\eps -(1/d)})
+ O (n^{(\alpha(d-2)/(2d(d-1))) - (1/2) +(2+6\alpha)\eps} )
,\]
which again tends to zero for $\eps$ small enough, giving (\ref{101}). $\square$\\

Recall the definition of the point process $\V_n \subset (0,1)^{d-1}$
from (\ref{vndef}).

\begin{lemma}
\label{kk1}
Suppose $d \in \{2,3,\ldots\}$.
For $\alpha \in (0,d)$ we have that as $n \to \infty$
\begin{align}
\label{968}
 n^{(\alpha/d)-1} \O^{d-1,\alpha} (\V_n) \to 0, \textrm{ in } L^1.\end{align}
Also, for $\alpha \geq d$, there exist finite
positive
constants $\mu(d-1,\alpha)$ such that as $n \to \infty$,
\begin{align}
\label{969}
 \Exp [ \O^{d-1,\alpha} (\V_n) ] \to \mu(d-1,\alpha).\end{align}
 In particular, $\mu(1,\alpha) = \frac{2}{\alpha(\alpha+1)} ( 1+\frac{2^{-\alpha}}{\alpha-1} )$.
\end{lemma}
\proof Suppose $\alpha \in (0,d)$. 
Recall that $\beta_n = \card (\V_n)$ is Poisson with mean $\Theta(n^{(1/2)-\eps})$.
Let $\bU_1,\bU_2,\ldots$ be a sequence of independent uniform random vectors on $(0,1)^d$.
 Let $\U_m$ denote the 
sequence
of uniform random vectors in $(0,1)^{d-1}$
formed by the sequence orthogonal projections down onto $(0,1)^{d-1}$ of the
points of $\{ \bU_1,\ldots,\bU_m\} \cap B_n$
listed in order of increasing $x_d$-coordinate. Then, without loss of generality, 
we can assume that 
$\Po_n = \{ \bU_1,\ldots, \bU_{N(n)}\}$
with $N(n)$ Poisson with mean $n$,  
$\beta_n = \card (\Po_n \cap B_n)$, and
$\V_n = \U_{\beta_n}$ in this notation.
 
Let $A_n$ denote the event
$\{\beta_n > n t_n +n^{1/4}\}$. Then by standard
Chernoff bounds on Poisson tails (see, e.g., Lemma 1.2 of \cite{penbook}),
$\Pr(A_n) = O( \re^{-Cn^\eps})$ for some $C \in (0,\infty)$. With the coupling described above,
\begin{align}
\label{pp1}
 n^{(\alpha/d)-1} \O^{d-1,\alpha} (\V_n) \leq n^{(\alpha/d)-1} \O^{d-1,\alpha} (\U_{ \lceil
 n t_n
 +n^{1/4} \rceil })
+n^{(\alpha/d)-1} \1_{A_n} C' N(n) ,\end{align}
for some $C' \in (0,\infty)$
and $N(n) = \card (\Po_n)$ is Poisson with mean $n$. By  Theorem 2.1 of \cite{ong},
for $\alpha < d-1$
we have that as $m \to \infty$
\begin{align}
\label{9636}
 \Exp [ \O^{d-1,\alpha} (\U_m) ] = O( m^{(d-1-\alpha)/(d-1)}),\end{align}
and also
\begin{align}
\label{967}
 \Exp[ \O^{d-1,d-1} (\U_m)] = O(\log m), ~~~ \Exp[ \O^{d-1,\alpha} (\U_m)] \to \mu(d-1,\alpha)~~~(\alpha>d-1),\end{align}
for some positive constant $\mu(d-1,\alpha)$:
this notation coincides
with Proposition 2.1 of \cite{ong2}. The particular
values $\mu(1,\alpha) = \frac{2}{\alpha(\alpha+1)} ( 1+\frac{2^{-\alpha}}{\alpha-1} )$ for
$\alpha >1$ were given in Proposition 2.1 of \cite{ong}. Thus by (\ref{9636}), if $\alpha < d-1$,
\[ \Exp [ n^{(\alpha/d)-1} \O^{d-1,\alpha} (\U_{
\lceil n t_n +n^{1/4} \rceil}) ]
= O \left(n^{-(1/2)-\eps +\alpha((d-2+\eps d)/(2d(d-1)))} \right) \to 0,\]
as $n\to \infty$, for $\eps$ small. Also, for $\alpha \in [d-1,d)$
\[ \Exp [ n^{(\alpha/d)-1} \O^{d-1,\alpha} (\U_{\lceil
n t_n +n^{1/4} \rceil}) ]
 \to 0,\]
 as $n \to \infty$, by (\ref{967}). 
 Also by Cauchy--Schwarz
 \begin{align}
 \label{pp2}
  \Exp [ n^{(\alpha/d)-1} \1_{A_n}  N(n) ] \leq n^{(\alpha/d)-1} (\Exp [ N(n)^2] )^{1/2}
 (\Pr (A_n))^{1/2} \to 0 ,\end{align}
 as $n \to \infty$. So from (\ref{pp1}) this completes the proof of (\ref{968}).
 
For the proof of (\ref{969}), let $A'_n$ denote the event that $\{\beta_n < n t_n -n^{1/4}\}$. Then
by Chernoff tail bounds again, $\Pr(A'_n) = O(\re^{-Cn^\eps})$.
We have that there is a constant $C' \in (0,\infty)$ such that for all $n$
\begin{align}
\label{965}
 \O^{d-1,\alpha} (\U_{\lceil n t_n -n^{1/4} \rceil}) -\1_{A'_n} C' n \leq
\O^{d-1,\alpha} (\V_n)
\leq
\O^{d-1,\alpha} (\U_{\lceil n t_n+n^{1/4} \rceil}) +\1_{A_n} C' N(n).\end{align}
Suppose $\alpha \geq d > d-1$. Then by (\ref{967}) and (\ref{pp2})
we have that
the expectations of both the lower and upper bounds in (\ref{965})
converge to $\mu(d-1,\alpha)$. Thus we have (\ref{969}). $\square$\\

\noindent
{\bf Proof of Theorem \ref{llnthm}.} Consider
\begin{align}
\label{000} 
\LL^{d,\alpha} (\Po_n) = \LL^{d,\alpha} (\Po_n ; \Gamma_n)
+ \LL^{d,\alpha} (\Po_n ; B_n)
+ \LL^{d,\alpha} (\Po_n ; I_n) .\end{align}
First suppose $\alpha \in (0,d)$. We have
\begin{align}
\label{8888}
 \Exp [ n^{(\alpha/d)-1} \LL^{d,\alpha} (\Po_n ; B_n) ] \nonumber\\
= \Exp [ n^{(\alpha/d)-1} \O^{d-1,\alpha} (\V_n) ]
+ n^{(\alpha/d)-1} \Exp \left[  \LL^{d,\alpha} (\Po_n ; B_n)
- \O^{d-1,\alpha} (\V_n) \right].\end{align}
From (\ref{968}) we have that the first term on the right-hand side
of (\ref{8888})
tends to zero as $n \to \infty$ for $\alpha \in (0,d)$. By
(\ref{0214aa}) we have that for $\alpha \in (0,1)$
the second term on the right-hand side of (\ref{8888})
is $O(n^{\alpha((1/d)-(1/2)-\eps)-(1/2)-\eps})$
which tends to zero for $d \geq 2$, and
(\ref{0209cc}) yields the same result for $\alpha \geq 1$. 
Thus for any $\alpha \in (0,d)$, we have that $n^{(\alpha/d)-1} \LL^{d,\alpha} (\Po_n ; B_n)$
tends to zero in $L^1$.
Then multiplying
both sides of (\ref{000}) by $n^{(\alpha/d)-1}$ and applying Lemma \ref{llnlem} and (\ref{101})
we obtain (\ref{0728e}).  

Now suppose $\alpha \geq d$. We have
\begin{align}
\label{8888a}
\Exp [ \LL^{d,\alpha} (\Po_n ; B_n) ]
= \Exp [ \O^{d-1,\alpha} (\V_n) ]
+ \Exp \left[  \LL^{d,\alpha} (\Po_n ; B_n)
- \O^{d-1,\alpha} (\V_n) \right].\end{align}
By (\ref{0209cc}) 
 the last term on the right
of (\ref{8888a}) tends to zero as $n \to \infty$, since $\alpha \geq d > 1$.
Also, (\ref{969}) says that the first term on the right of (\ref{8888a})
tends to $\mu(d-1,\alpha)$. Thus for $\alpha \geq d$
\[ \Exp [ \LL^{d,\alpha} (\Po_n ; B_n) ] \to \mu(d-1,\alpha).\]
Also, Lemma \ref{llnlem} implies that, for $\alpha > d$
\[ \Exp [ \LL^{d,\alpha} (\Po_n ; \Gamma_n ) ] \to 0,\]
and for $\alpha =d$
\[ \Exp [ \LL^{d,d} (\Po_n ; \Gamma_n) ] \to 2 v_d^{-1} .\]
Then taking expectations in (\ref{000})
and using (\ref{1111d}) gives (\ref{011b}). This completes the proof
of Theorem \ref{llnthm}. $\square$

\end{document}